\documentclass[%
 aip,
 amsmath,amssymb,
 reprint,%
]{revtex4-2}

\usepackage{graphicx}
\usepackage{dcolumn}
\usepackage{bm}

\usepackage[utf8]{inputenc}
\usepackage[T1]{fontenc}
\usepackage{mathptmx}
\usepackage{etoolbox}
\usepackage{graphicx}
\usepackage{dcolumn}
\usepackage{bm}
\usepackage{enumitem}
\usepackage{amssymb}
\usepackage{amsmath}
\usepackage{psfrag}
\usepackage{mathtools}
\usepackage{color}
\usepackage{xcolor}
\usepackage[export]{adjustbox}
\usepackage[english]{babel}
\newtheorem{theorem}{Theorem}
\usepackage{hyperref}
\usepackage{natbib}
\usepackage{stmaryrd}

\definecolor{cyan1}{rgb}{0.387, 0.82, 1}
\definecolor{red1}{rgb}{0.902, 0.383, 0.355}
\definecolor{RGBred}{rgb}{1,0,0}
\definecolor{RGBblue}{rgb}{0,0,1}
\definecolor{outgreen}{rgb}{0.38,0.73,0.66}
\definecolor{arrowgreen}{rgb}{0,0.631,0.294}
\definecolor{arrowblue}{rgb}{0.129,0.251,0.6}
\definecolor{arrowgreen}{rgb}{0,0.631,0.294}
\definecolor{arrowred}{rgb}{0.929,0.11,0.141}
\definecolor{poleorange}{rgb}{0.984,0.69,0.231}
\definecolor{poleblue}{rgb}{0.161,0.667,0.886}
\definecolor{planegray}{rgb}{0.502,0.502,0.502}
\definecolor{orbitorange}{cmyk}{0,0.5,1,0}
\definecolor{bblue}{cmyk}{1,1,0,0}

\newcommand{\dc}[1]{{#1}}
\newcommand{\lc}[1]{{#1}}
\newcommand{\cc}[1]{{#1}}
\newcommand{\qq}[1]{{#1}}
\newcommand{\mm}[1]{{#1}}
\newcommand{\rr}[1]{{#1}}

\newtheorem{corollary}{Corollary}

\renewcommand{\arraystretch}{1.3}

\makeatletter
\def\@email#1#2{%
 \endgroup
 \patchcmd{\titleblock@produce}
  {\frontmatter@RRAPformat}
  {\frontmatter@RRAPformat{\produce@RRAP{*#1\href{mailto:#2}{#2}}}\frontmatter@RRAPformat}
  {}{}
}%
\makeatother
\begin{document}

\preprint{AIP/123-QED}

\title{\mm{Smooth transformations and ruling out closed orbits in planar systems}}

\author{Tiemo Pedergnana}
 \email{ptiemo@ethz.ch}
 \affiliation{%
 CAPS Laboratory, Department of Mechanical and Process Engineering, ETH Z{\"u}rich, Sonneggstrasse 3,
8092 Z{\"u}rich, Switzerland
}%
\author{Nicolas Noiray}%
 \email{noirayn@ethz.ch}
\affiliation{%
 CAPS Laboratory, Department of Mechanical and Process Engineering, ETH Z{\"u}rich, Sonneggstrasse 3,
8092 Z{\"u}rich, Switzerland
}%

\date{\today}

\begin{abstract}
\mm{This work deals with planar dynamical systems with and without noise. In the first part, we seek to gain a refined understanding of such systems by studying their differential-geometric transformation properties under an arbitrary smooth mapping. Using elementary techniques, we obtain a unified picture of different classes of dynamical systems, some of which are classically viewed as distinct. We specifically give two examples of Hamiltonian systems with first integrals, which are simultaneously gradient systems. Potential applications of this apparent duality are discussed. The second part of this study is concerned with ruling out closed orbits in steady planar systems. We reformulate Bendixson's criterion using the coordinate-independent Helmholtz decomposition derived in the first part, and we derive another, similar criterion. Our results allow for automated ruling out of closed orbits in certain regions of phase space, and could be used in the future for efficient seeding of initial conditions in numerical algorithms to detect periodic solutions.}
\end{abstract}

\maketitle

\begin{quotation}
\mm{Planar dynamical systems can exhibit a broad range of complicated dynamics, some of which remain elusive even today. While great simplifications are available for Hamiltonian or gradient systems, the general Helmholtz decomposition has so far found little application in nonlinear dynamics. Since dynamical systems can be defined ad-hoc, for example, to model phase transitions or biological processes, the choice of basis for a given system may be ambiguous, which makes identification of the Helmholtz decomposition a nontrivial task. In contrast, mechanical systems generally have a known, preferred basis. In the first part of this work, we study bivariate Langevin equations under an arbitrary smooth mapping to obtain the transformation properties of their Helmholtz decomposition. This investigation reveals a unified picture of different classes of planar systems, some of which are are typically presented as distinct. Specifically, we give two explicit examples of Hamiltonian systems with first integrals which are simultaneously gradient systems. In the second part, we discuss criteria for ruling out closed orbits in steady planar systems. We reformulate Bendixson's criterion in terms of the coordinate-independent Helmholtz decomposition derived in the first part, and we present another criterion for ruling out closed orbits. These latter results could help simplify numerical algorithms to detect periodic solutions in planar flows.}
\end{quotation}

\section{Overview}
\subsection{Dynamical system}
In this work, we study planar (stochastic) dynamical systems given by a bivariate Langevin equation of the form\cite{stratonovich1963topics,gardiner1985handbook,risken1996fokker}
\begin{eqnarray}
    \dot{x}=\mathcal{F}(x,t)+\mathcal{B}(x)\Xi, \label{Dynamical system}
\end{eqnarray}
which is defined on a subset of the plane $\mathcal{D}\subset \mathbb{R}^2$.\footnote{ Throughout this work, a dot over a dependent variable denotes the total time derivative.} Equation \eqref{Dynamical system} states that the evolution of the random variable $x$: $\mathbb{R}^+\rightarrow \mathbb{R}^{2}$ over time $t\in\mathbb{R}^+$ is governed by the smooth velocity field $\mathcal{F}(x,t)$: $\mathbb{R}^2\times\mathbb{R}^+\rightarrow \mathbb{R}^{2}$, the diffusion tensor $\mathcal{B}$: $\mathbb{R}^2\times\mathbb{R}^+\rightarrow \mathbb{R}^{2\times 2}$ and the vector $\Xi=(\xi_1,\xi_2)^T$, whose entries are white Gaussian noise sources with variance $\Gamma$ and zero mean.\cite{stratonovich1963topics} Note that, if either $\Gamma=0$ or $\mathcal{B}$ vanishes identically for all $x\in\mathcal{D}$, the system given by Eq. \eqref{Dynamical system} defines a deterministic planar flow with velocity field $\mathcal{F}$.\cite[see pp. 42--65 and pp. 125--306, respectively]{guckenheimer2013nonlinear,strogatz2018nonlinear} \qq{Unless explicitly stated otherwise,} all quantities are assumed to be real in this work. 

\subsection{Helmholtz decomposition}
 \mm{Helmholtz's theorem\mm{\cite{stokes_2009,Helmholtz185825}} states that, under certain conditions,} \qq{a} smooth, planar velocity field $u$ can be decomposed as follows:\cite[see][pp. 52--54]{morse1953methods}
\begin{eqnarray}
u(x,t)&=&{-\nabla \mathcal{V}(x,t)}+{S\nabla \mathcal{H}(x,t)}, \label{Helmholtz decomposition}
\end{eqnarray}
where $[\nabla (\cdot)]_m={\partial (\cdot)}/{\partial{x_m}}$, $\mathcal{V}$  is the scalar potential, $\mathcal{H}$ is the Hamiltonian function and 
\begin{eqnarray}
    S&=&\begin{pmatrix}
        0 & 1\\
        -1 & 0
    \end{pmatrix}. \label{Skew symm matrix}
\end{eqnarray}
The two-dimensional (2D) Helmholtz decomposition (HD) \eqref{Helmholtz decomposition} can be derived from the three-dimensional case given in the above references by setting the vector potential equal to $(0,0,\mathcal{H})^T$. \mm{We call $\mathcal{H}$ the Hamiltonian function because for $\mathcal{V}=0$, Eq. \eqref{Helmholtz decomposition} corresponds to the definition of a canonical Hamiltonian velocity field.\cite{arnol2013mathematical} The restriction to the 2D case is motivated by the fact that for planar vector fields, the Helmholtz decomposition requires no special gauge conditions, since the number of potentials is equal to the number of entries of $u$.} \mm{Nevertheless, there is some ambiguity in the decomposition \eqref{Helmholtz decomposition}}: Adding a constant to $\mathcal{V}$ or $\mathcal{H}$ \qq{does not change $u$ defined in Eq. \eqref{Helmholtz decomposition}.} Furthermore, a \mm{spatially} constant term on the right-hand-side (RHS) of Eq. \eqref{Dynamical system} can be included as a linear (monomial) term either in $\mathcal{V}$ or in $\mathcal{H}$. See Appendix \ref{Example Harmonic Oscillator} for an example. \mm{No such ambiguity occurs in higher-order polynomial terms in a broad class of planar systems representing nonlinear oscillators, as shown in Appendix \ref{Appendix C}.} The HD is broadly used in fluid mechanics, \cite{Morino198665,Joseph200614272,Linke2014782,Buhler20141007,Lindborg2014,Buhler2017361,Schoder20203019,Caltagirone2021} electromagnetics, \cite{Aharonov1959485,Konopinski1978499,Haber2000150,hehl2003foundations} geophysics, \cite{Weiss201340,Du2017S111,Shi2019509} imaging \cite{Paganin19982586,Park20063697} and computer vision. \cite{Kohlberger2003432,Guo2005493,Cuzol2007329} The reader interested in mathematical discussions of the HD, \cite{simader92,Farwig2007239,sohr2012navier,moses1971,Schweizer2018,Giga2022} its generalization to $n$-dimensional vector fields\cite{Glotzl2023} or a review of its applications\cite{Bhatia20131386} is referred to the respective literature. Uniqueness or existence are not of concern in this work, which is focused on the transformation properties of the HD  \eqref{Helmholtz decomposition} for smooth planar fields. Therefore, in the following, we generally assume that a (quasi-)unique decomposition exists throughout $\mathcal{D}$, such that Eq. \eqref{Helmholtz decomposition} serves as a definition of $u$. 


\subsection{Transformation properties \label{Intro part 3}}
\begin{figure}[t]
\begin{psfrags}
\psfrag{a}{$x=f(y)$}
\psfrag{b}{\hspace{-0.1cm}$\mathcal{D}\hspace{3.7cm}f^{-1}(\mathcal{D})$}
\psfrag{d}{}
\psfrag{c}{\hspace{0.15cm}\textcolor{red}{$C$}}%
\psfrag{d}{\textcolor{bblue}{$dl$}}
\begin{center}
\includegraphics[width=0.3\textwidth]{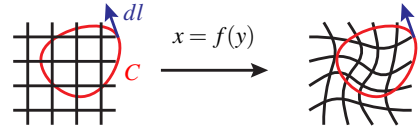}
\end{center}
\end{psfrags}
\vspace{-0.5cm}
    \caption{\mm{This work is concerned with smooth transformations and ruling out closed orbits in planar dynamical systems $\dot{x}=u(x,t)$. An example of a closed orbit is shown in red in the figure, superimposed on a transformed coordinate grid (black). Shown in blue is the tangent vector $dl$ to the curve $C$ at a certain point.}}
    \label{Figure 1}
\end{figure}

In expressing the HD by Eq. \eqref{Helmholtz decomposition}, the use of Cartesian (rectangular) coordinates was tacitly assumed. \cite[see pp. 21--52]{morse1953methods} With this in mind, we can study the dynamics of a planar dynamical system under additive white noise, $\dot{x}=u+\Xi$, in such coordinates. By Eq. \eqref{Helmholtz decomposition}, such a system can be written as
\begin{eqnarray}
\dot{x}&=&{-\nabla \mathcal{V}(x,t)}+{S\nabla \mathcal{H}(x,t)}+\Xi. \label{Cartesian system}
\end{eqnarray}
At first sight, comparing Eqs. \eqref{Dynamical system} and \eqref{Cartesian system} is sensible only in the special case when $\mathcal{B}$ corresponds to the 2-by-2 identity matrix. However, under an arbitrary, smooth mapping (see Fig. \ref{Figure 1})
\begin{eqnarray}
    x=f(y) \label{General mapping},
\end{eqnarray}
after redefining $y\rightarrow x$, Eq. \eqref{Cartesian system} is transformed into the following system:
\begin{eqnarray}
    \dot{x}=-g^{-1}(x)\nabla \widetilde{\mathcal{V}}(x,t) \pm\dfrac{1}{\sqrt{\det[g(x)]}}S\nabla\widetilde{\mathcal{H}}(x,t)+h(x)^{-1}\widetilde{\Xi}, \nonumber \\
    \label{Transformed 2D Helmholtz decomposition}
\end{eqnarray}
where $J=\nabla f$ is the Jacobian matrix of the mapping $f$, $J=Qh$ is the polar decomposition of $J$,\cite[see p. 449]{horn2012matrix} $Q=Q^{-T}\in\mathbb{R}^{2\times 2}$ is orthogonal, $h\in\mathbb{R}^{2\times 2}$ is positive definite, $(\cdot)^{T}$ is the transpose, $g=h^T h\in\mathbb{R}^{2\times 2}$ is the positive definite metric tensor of the mapping $f$, $\det(\cdot)$ is the determinant, $\widetilde{\mathcal{V}}(x,t)=\mathcal{V}(f(x),t)$ is the transformed scalar potential and $\widetilde{\mathcal{H}}(x,t)=\mathcal{H}(f(x),t)$ is the transformed Hamiltonian. The ``$\pm$''-sign in Eq. \eqref{Transformed 2D Helmholtz decomposition} indicates whether $Q$ is purely rotational ($\det Q=1$, ``$+$'') or if it contains a reflection ($\det Q=-1$, ``$-$''). By the positive definiteness of $h$, $\det h(x)=\sqrt{\det g(x)}$ is positive. \footnote{Note that, by the uniqueness of the Cholesky decomposition for positive definite matrices,\cite[see p. 441]{horn2012matrix} the polar decomposition $J=Qh$ and the QR factorization\cite[see p. 449]{horn2012matrix} of $J$ coincide, which implies that $h^{-1}$ is generally a lower triangular matrix.} We also mention the implicit assumption in the derivation of Eq. \eqref{Transformed 2D Helmholtz decomposition} that the transformed noise term $\widetilde{\Xi}$ is related to its original counterpart by\cite[see Sec. III. B.]{Pedergnana2022} $\Xi=Q(x)\widetilde{\Xi}$. \mm{It is understood that the expression on the RHS of Eq. \eqref{Transformed 2D Helmholtz decomposition} is valid only in an open neighborhood in which $g$ is non-singular. See pp. 172--182 of \citet{guillemin2010differential} for a detailed discussion. The transformation formula defined by Eq. \eqref{Transformed 2D Helmholtz decomposition} is exemplified on the noise-driven, forced-damped Harmonic oscillator in Sec. \ref{Appendix B}.}

The first and last terms on the RHS of Eq. \qq{\eqref{Transformed 2D Helmholtz decomposition}} have been derived for arbitrary-dimensional systems in previous work.\cite{Pedergnana2022} The present study complements those results for the planar case by adding the transformed Hamiltonian term. If we drop all other terms and assume $\det Q=1$, the transformed Hamiltonian system reads
\begin{eqnarray}
    \dot{x}=\dfrac{1}{\sqrt{\det[g(x)]}}S\nabla\widetilde{\mathcal{H}}(x,t). \label{Transformed Hamiltonian system}
\end{eqnarray}
which is consistent with earlier works.\cite{Gonzalez-Gascon198661,Nutku199027} \mm{If a given system can be identified to be of the form \eqref{Transformed Hamiltonian system}, there exist sets of preferred coordinates for that system, namely those in which $\det g$ becomes unity, leading to a canonical Hamiltonian system.}


\subsection{Steady planar flows}
\mm{An important feature of steady planar flows
\begin{eqnarray}
    \dot{x}=u(x), \label{Steady planar system}
\end{eqnarray}
is the existence or absence of a closed orbit $C$ in the phase space (see Fig. \ref{Figure 1}).\footnote{An ``orbit'' is \mm{a curve} in the (frozen) phase space of a steady dynamical system along which a trajectory runs.} Determining the number of limit cycles in a given steady planar system is a nontrivial task\cite{delpino2023limit} which can be tackled using numerical methods \cite{parker2012practical} or bifurcation theory\cite{han2006bifurcation,guckenheimer2013nonlinear}. Determining an upper limit to the number of limit cycles in polynomial systems of the form \eqref{Steady planar system} is also an open question and part of Hilbert's 16\textsuperscript{th} problem.\cite{hilbert1900mathematische,Ilyashenko2002301}}

\mm{There exist few simple, analytical criteria to rule out closed orbits in steady planar systems.\cite[see p. 44 and pp. 201--205, respectively]{guckenheimer2013nonlinear,strogatz2018nonlinear} One example of such a criterion is Dulac's criterion, which states that if there exists a smooth function $\varphi=\varphi(x)$ (the Dulac function\cite{Busenberg1993463}) such that the divergence of $\varphi u$ has only one sign throughout a simply connected domain $\mathcal{D}$, then Eq. \eqref{Steady planar system} has no closed orbits\footnote{The \mm{limiting} case of a closed orbit which is a fixed point is excluded from the discussion here.} which are entirely contained $\mathcal{D}$.\cite[see][p. 204]{strogatz2018nonlinear} In the special case where $\varphi=\mathrm{const.}$, Dulac's criterion reduces to Bendixson's criterion.\cite[see p. 44]{guckenheimer2013nonlinear}}

\mm{Combining the main result from the previous section, Eq. \eqref{Transformed 2D Helmholtz decomposition}, with the theory of \citet{jost2008riemannian}, we show in Appendix \ref{Appendix A} that Bendixson's criterion can be reformulated as a statement which only concerns the scalar potential $\widetilde{\mathcal{V}}$ of $u$:}

\mm{\begin{theorem}[Reformulation of Bendixson's criterion]\label{Theorem 2}
If the scalar potential $\widetilde{\mathcal{V}}$, or its negative $-\widetilde{\mathcal{V}}$, of the planar field 
\begin{eqnarray}
    u(x)=-g^{-1}(x)\nabla \widetilde{\mathcal{V}}(x) \pm\dfrac{1}{\sqrt{\det[g(x)]}}S\nabla\widetilde{\mathcal{H}}(x), \label{steady planar field}
\end{eqnarray}
is a strictly subharmonic function throughout a simply connected subset of the plane $\mathcal{D}\subset \mathbb{R}^2$, then the system $\dot{x}=u(x)$, $u=(u_1,u_2)^T$, has no closed orbits contained entirely in $\mathcal{D}$.
\end{theorem}
This result suggests that the Helmholtz decomposition is, to some degree, meaningful for (steady) planar systems, as the function $\widetilde{\mathcal{V}}$ alone can rule out the existence of closed orbits, regardless of the form of $\widetilde{\mathcal{H}}$.}

\mm{In Sec. \ref{Main results 2}, we analytically derive another criterion for ruling out closed orbits. This criterion is first tested on classic examples of linear and nonlinear oscillators and then on the system of Shi Singlong,\cite{shi1980concrete,Shi2019509} an example of a quadratic planar system with exactly four limit cycles.\cite{Galias2022} We note that four is the maximal number of limit cycles in a quadratic planar system known to date.\cite{Ilyashenko2002301}}

\section{Relation to prior work \label{prior work section}}
\mm{There exists a substantial literature on decompositions of smooth vector fields and on dynamical systems with Hamiltonian, gradient or mixed structure. This section seeks to connect the findings of the present work to those prior efforts.}
\mm{The derivation of the transformation formula in Eq. \eqref{Transformed 2D Helmholtz decomposition} partly overlaps with Sec. 3.4.2 of \citet{risken1996fokker}. The present work differs from Risken's derivation in that, by writing out the Helmholtz decomposition explicitly, further simplifications are enabled. For example, the orthogonal matrix $Q$ from the polar decomposition of the Jacobian essentially drops out of the transformed dynamics, which is not evident in the reference. In this regard, we also mention the Lamperti transform \cite[see][p. 98--100]{sarkka2019applied}, which can be used to transform multiplicative to additive noise in one-dimensional stochastic systems.}

\mm{The coordinate-independent formulation of the Helmholtz decomposition in Eq. \eqref{Transformed 2D Helmholtz decomposition} (set $\widetilde{\Xi}=0$ for comparison) is a special case of the Helmholtz--Hodge decomposition theorem \cite[see][p. 539]{abraham2012manifolds}, applied to planar systems. Other concepts related to this work are Poisson structures and Poisson brackets\cite[see][pp. 381--387 and pp. 456--468]{olver1993applications,arnol2013mathematical}. In the language of \citet{olver1993applications}, for example, the matrix $\sqrt{\det g}^{-1}S$ is the structure matrix for the Poisson bracket associated with the Hamiltonian system described by Eq. \eqref{Transformed Hamiltonian system} in the steady case. \citet{Quispel1996223} denote (steady) systems of the form of Eq. \eqref{Transformed Hamiltonian system} simply as ``Poisson systems''. \citet{McLachlan19982399} show that any system with a first integral (a conserved quantity) can be written as a Poisson system (see \citet{McLachlan19991021} for the corresponding proofs). By combining these prior results with the present work, it becomes evident that, in the planar case, systems with first integrals are mapped into canonical, steady Hamiltonian systems by smooth transformations.}

\mm{\citet{McLachlan19982399} further show that any system with a strict (strictly decreasing along trajectories) Lyapunov function can be written as $\dot{x}=-g^{-1}(x)\nabla \widetilde{V}$, where $g$ is a positive definite matrix. Such systems have been shown to derive from gradient-driven dynamics, subjected to an arbitrary smooth mapping $f$.\cite{Pedergnana2022} The combination of the above results implies that smooth transformations map systems with strict Lyapunov functions into canonical gradient systems, which is also consistent with the work of \citet{Barta201257}. }

\begin{figure}[t]
\begin{psfrags}
\psfrag{a}{\hspace{-0.03cm}$\beta=-1$\hspace{1.5cm}$\beta=1$}
\psfrag{b}{$\mathcal{H}$}
\psfrag{c}{\textcolor{red}{$\widetilde{\mathcal{V}}$}}

\begin{center}
\includegraphics[width=0.25\textwidth]{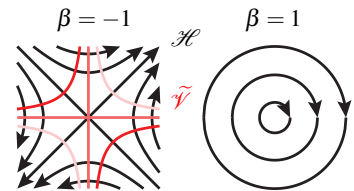}
\end{center}
\end{psfrags}
\vspace{-0.5cm}
    \caption{\mm{For $\beta>0$, Eq. \eqref{HO Ham. sys} describes the Harmonic oscillator. For $\beta<0$, this system can simultaneously be described as governed by a Hamiltonian $\mathcal{H}$ or as driven by the gradient of a scalar potential $\widetilde{\mathcal{V}}$. Superimposing the isocontours of $\mathcal{H}$ (black) and $\widetilde{\mathcal{V}}$ (red, more saturated colors correspond to larger values) for $\beta=-1$ reveals an intuitive picture: the potential landscape drives trajectories along the isocontours of the Hamiltonian. Arrows indicate the direction of the flow.}}
    \label{Figure 2}
\end{figure}

\mm{In standard textbooks on nonlinear dynamics, gradient systems and Hamiltonian systems are presented separately.\cite{guckenheimer2013nonlinear,strogatz2018nonlinear} However, assuming these classes of systems to be inherently distinct would be a false dichotomy. In general, a given system can \textit{simultaneously} be a gradient system and a Hamiltonian system, which is already evident in the trivial system $\dot{x}=0$. For an explicit, nontrivial example, consider the canonical Hamiltonian system \begin{eqnarray}
    \dot{x}=\begin{pmatrix}
        0 & 1\\
        -1 & 0
    \end{pmatrix}\nabla \mathcal{H}(x), \label{HO Ham. sys}
\end{eqnarray}
$x=(x_1,x_2)^T$, with $\mathcal{H}(x)=\beta x_1^2/2+x_2^2/2$. If $\beta>0$, this system corresponds to the Harmonic oscillator. However, for \textit{any} $\beta$, this system is Hamiltonian and has a first integral: $d\mathcal{H}/dt=0$. Now, note that the same system can be rewritten as
\begin{eqnarray}
    \dot{x}=-\underbrace{\begin{pmatrix}
        1 & 0\\
        0 & -\beta
    \end{pmatrix}}_{=g^{-1}(x)}\nabla \widetilde{\mathcal{V}}(x), \label{HO gradient sys}
\end{eqnarray}
where $\widetilde{\mathcal{V}}=-x_1 x_2$. For $\beta<0$, the matrix $g^{-1}=\mathrm{diag}(1,-\beta)$ is positive definite, i.e., Eq. \eqref{HO gradient sys} represents a (transformed) gradient system.\cite{Pedergnana2022} Note that, in this parameter range, the isocontours of $\mathcal{H}$ describe hyperbolas in the plane, as shown in Fig. \ref{Figure 2}. For $\beta=-1$, $g$ becomes equal to the identity matrix, and the contours of $\widetilde{\mathcal{V}}$ are exactly orthogonal to those of $\mathcal{H}$: the potential landscape drives the system's trajectories along the energy levels of the Hamiltonian. The same is not true for $\beta\neq -1$, as $g$ skews the gradient dynamics, leading the trajectories to not follow exactly the steepest descent. For $\beta<0$, the real mapping $f(x)=(x_1,x_2 \sqrt{-\beta})^T$, for instance, transforms Eq. \eqref{HO gradient sys} back into a canonical gradient system $\dot{x}=-\nabla \mathcal{V}$, and Eq. \eqref{HO Ham. sys} into $\dot{x}=\sqrt{\det g}^{-1} S\nabla \widetilde{\mathcal{H}}$. In these preferred coordinates, the isocontours of ${\mathcal{V}}$ are always perpendicular to those of $\widetilde{\mathcal{H}}$. Another example of a system which is simultaneously a gradient system and a Hamiltonian system is given in Sec. \ref{Swarm system}.}


\section{Main results \label{Main results}}
\subsection{Transformed Hamiltonian systems \label{Main results 1}}
Since the gradient and noise terms in Eq. \eqref{Transformed 2D Helmholtz decomposition} have already been derived,\cite{Pedergnana2022} what is left is to obtain the transformation formula for the Hamiltonian system $\dot{x}=S\nabla \mathcal{H}$ under a general mapping of the form \eqref{General mapping}. For this, we set $x=f(y)$ and apply the chain rule to this system in index form:
\begin{eqnarray}
    \dot{x}_m&=&S_{mn}\frac{\partial \mathcal{H}(f(y),t)}{\partial x_n},\\
    \frac{\partial x_m}{\partial y_k} \dot{y}_k &=&S_{mn}\frac{\partial y_l}{\partial x_n}\frac{\partial \widetilde{\mathcal{H}}(y,t)}{\partial y_l}, \\
   \dot{y}_k &=&\underbrace{\frac{\partial y_k}{\partial x_m}S_{mn}\frac{\partial y_l}{\partial x_n}}_{\widetilde{S}_{kl}}\frac{\partial \widetilde{\mathcal{H}}(y,t)}{\partial y_l}, \label{int res}
\end{eqnarray}
where the time-dependence of $x$ and $y$ was suppressed for brevity. In Eq. \eqref{int res}, we defined the matrix $\widetilde{S}$, whose entries are given by
\begin{eqnarray}
    \widetilde{S}_{kl}&=&\frac{\partial y_k}{\partial x_m}S_{mn}\frac{\partial y_l}{\partial x_n}\\
    &=&\frac{\partial y_k}{\partial x_1}\frac{\partial y_l}{\partial x_2}-\frac{\partial y_k}{\partial x_2}\frac{\partial y_l}{\partial x_1}\\
    \implies \widetilde{S}&=&\Big(\frac{\partial y_1}{\partial x_1}\frac{\partial y_2}{\partial x_2}-\frac{\partial y_1}{\partial x_2}\frac{\partial y_2}{\partial x_1}\Big)S. \label{int res 2}
\end{eqnarray}
To interpret this result, we note that the term in brackets is simply the determinant of the Jacobian of the inverse mapping $y=f^{-1}(x)$. By the properties of the determinant,\cite[see pp. 8--12]{horn2012matrix} we conclude that 
\begin{eqnarray}
    \frac{\partial y_1}{\partial x_1}\frac{\partial y_2}{\partial x_2}-\frac{\partial y_1}{\partial x_2}\frac{\partial y_2}{\partial x_1}=[\det J(y)]^{-1}.
\end{eqnarray}
A further simplification is enabled by the polar decomposition $J=Qh$, where $Q$ is an orthogonal matrix with $\det(Q)=\pm 1$. Note that $\det{J}=\det{Q}\det{h}$ and $\det{h}=\sqrt{\det{g}}$. Collecting the above results, we have shown that under a general mapping $f$, a planar Hamiltonian system transforms like
\begin{eqnarray}
    S\nabla \mathcal{H}(x,t)\rightarrow \pm \dfrac{1}{\sqrt{\det[g(y)]}} S{\nabla}_y\,  \widetilde{\mathcal{H}}(y,t), \label{int res 3}
\end{eqnarray}
where \qq{$({\nabla}_y)_m=\partial /\partial y_m$} and the transformed Hamiltonian is defined as $\widetilde{\mathcal{H}}(y,t)=\mathcal{H}(f(y),t)$. Redefining $y\rightarrow x$, \qq{$\nabla_y\rightarrow \nabla$} in Eq. \eqref{int res 3} and combining this formula with previous results\cite{Pedergnana2022} yields Eq. \eqref{Transformed Hamiltonian system}.

\subsection{Ruling out closed orbits \label{Main results 2}}
Dulac's criterion concerns the autonomous planar dynamical system defined by Eq. \eqref{Steady planar system}, which is governed by the smooth vector field $u$ defined on a simply connected domain $\mathcal{D}\subset\mathbb{R}^2$. The criterion states that if there exists a smooth function $\varphi=\varphi(x)$ (the Dulac function\cite{Busenberg1993463}) such that the divergence of $\varphi u$ has only one sign throughout $\mathcal{D}$, then Eq. \eqref{Steady planar system} has no closed orbits which are entirely contained in that region.\cite[see][p. 204]{strogatz2018nonlinear} Here, under the same assumptions on $u$ and $\mathcal{D}$, we prove the following statement:

\begin{theorem}\label{Theorem 1}
    If there exists a smooth, positive definite $2$-by-$2$ matrix function $\mathcal{N}=\mathcal{N}(x)$ \qq{with entries $\mathcal{N}_{mn}(x)$, $m,n\in\{1,2\}$,} such that the \qq{out-of-plane} \mm{component} $\omega$ of
    \begin{eqnarray}
        \mathrm{curl}\,U(x)=\begin{pmatrix}
            0\\
            0\\
            \omega(x)
        \end{pmatrix}, \label{theorem 1 definition}
    \end{eqnarray}
    where
     \begin{eqnarray}
       U(x)=\begin{pmatrix}
            \mathcal{N}_{11}(x)u_1(x)+\mathcal{N}_{12}(x)u_2(x)\\
             \mathcal{N}_{21}(x)u_1(x)+\mathcal{N}_{22}(x)u_2(x)\\
            0
        \end{pmatrix}, \label{U definition}
    \end{eqnarray}is zero throughout a simply connected subset of the plane $\mathcal{D}\subset \mathbb{R}^2$, then the planar system $\dot{x}=u(x)$, $u=(u_1,u_2)^T$, has no closed orbits contained entirely in $\mathcal{D}$.
\end{theorem}
To prove Theorem \ref{Theorem 1}, we assume that the \qq{out-of-plane}-component of $\mathrm{curl}\,  U$, where $U$ is defined in Eq. \eqref{U definition}, is zero in $\mathcal{D}$. Integrating this expression over a simply connected subset $\mathcal{R}$ of the domain $\mathcal{D}$ bounded by the closed, positively oriented curve $C=\partial \mathcal{R}$, and applying Stokes' theorem\cite[see][p. 43]{morse1953methods} yields
\begin{eqnarray}
    0&=&\int_\mathcal{R} \omega(x) dA \label{criterion theorem 1}\\
    &=&\int_\mathcal{R} [\mathrm{curl}\, U(x)]^T n dA\\
    &=&\oint_C u^T(x) \mathcal{N}^T(x)  dl, \label{Dulac int res}
\end{eqnarray}
where $dA$ is an infinitesimal area element, $n=(0,0,1)^T$ is the normal vector and $dl$ is an infinitesimal \qq{tangent vector to $C$ (compare Fig. \ref{Figure 1}).} Along trajectories $x(t)$ governed by Eq. \eqref{Steady planar system}, $dl=udt$.\cite[see][p. 204]{strogatz2018nonlinear} Therefore, if $C$ is an orbit of the \qq{system given by Eq. \eqref{Steady planar system}, then it can be parametrized such that the RHS of Eq. \eqref{Dulac int res} takes the following form:}
\begin{eqnarray}
   \oint_C u^T(x(t)) \mathcal{N}^T(x(t))  u(x(t))dt. \label{Dulac int res 2}
\end{eqnarray}
By the positive definiteness of $\mathcal{N}$ ($	\Leftrightarrow$ of $\mathcal{N}^T$) and because $C$ is chosen arbitrarily, the integral \eqref{Dulac int res 2} generally has a positive sign,\footnote{Unless $u$ vanishes identically, in which case there is nothing to prove.} in contradiction to the \qq{assumption} that $\omega$ vanishes in $\mathcal{D}$. This contradiction shows that, under the assumptions of Theorem \ref{Theorem 1}, $\mathcal{D}$ can contain no closed curve $\Gamma$ which is an orbit of the system governed by Eq. \eqref{Steady planar system}. 

Theorem \ref{Theorem 1} can be applied in automated fashion by making, for example, a positive definite ansatz of the form
\begin{eqnarray}
    \mathcal{N}(x)=\begin{pmatrix}
        a(x_1) &\mathcal{N}_{12}(x_1,x_2)\\
        0& b
    \end{pmatrix}, \label{Ansatz 1}
    \end{eqnarray}
    or
    \begin{eqnarray}
    \mathcal{N}(x)=\begin{pmatrix}
        c &0\\
        \mathcal{N}_{21}(x_1,x_2)& d(x_2)
    \end{pmatrix}, \label{Ansatz 2}
\end{eqnarray}
where $a$ and $d$ are arbitrary, positive functions of $x_1$ and $x_2$, respectively, and $b$, $c\in\mathbb{R}^+$ are positive real numbers. We comment on this ansatz in Sec. \ref{Discussion section}. Note that $a$, $b$, $c$ and $d$ can be arbitrarily specified, and are only restricted by their positivity. In the examples we consider below, substituting Eqs. \eqref{Ansatz 1} and \eqref{Ansatz 2} into $U$ \qq{defined in Eq. \eqref{U definition}} and requiring that $\mathrm{curl}\, U=0$ yields two linear first-order differential equations for each system, one in $x_2$ (for $\mathcal{N}_{12}$) and one in $x_1$ (for $\mathcal{N}_{21}$). In each of these equations, the other respective coordinate is treated as a parameter. If the solution of either of the resulting equations exists over a simply connected planar domain $\mathcal{D}$, then, by Theorem \ref{Theorem 1}, there exists no closed orbit fully contained in $\mathcal{D}$. We note that, in a Cartesian basis, $\mathrm{curl}\,(\cdot)=\nabla \times (\cdot)$ and $\mm{\omega}=\partial U_1/\partial x_2-\partial U_2/\partial x_1$, so that the differential equations resulting from Eqs. \eqref{Ansatz 1} and \eqref{Ansatz 2} are
    \begin{eqnarray}
        -a(x_1)\dfrac{\partial u_1(x_1,x_2)}{\partial x_2}-\dfrac{\partial \mathcal{N}_{12}(x_1,x_2)}{\partial x_2}u_2(x_1,x_2)\nonumber\\
        -\mathcal{N}_{12}(x_1,x_2)\dfrac{\partial u_2(x_1,x_2)}{\partial x_2}+b\dfrac{\partial u_2(x_1,x_2)}{\partial x_1}=0, \label{corrolary 1}\\
        -c\dfrac{\partial u_1(x_1,x_2)}{\partial x_2}+\dfrac{\partial \mathcal{N}_{21}(x_1,x_2)}{\partial x_1}u_1(x_1,x_2)\nonumber\\
        +\mathcal{N}_{21}(x_1,x_2)\dfrac{\partial u_1(x_1,x_2)}{\partial x_1}+d(x_2)\dfrac{\partial u_2(x_1,x_2)}{\partial x_1}=0, \label{corrolary 2}
    \end{eqnarray}
respectively. The above results are summarized in the following corollary to Theorem \ref{Theorem 1}:
\begin{corollary} \label{Corollary 1}
    If the steady planar system $\dot{x}=u(x)$ is formulated in Cartesian coordinates, it has no closed orbits contained entirely in any simply connected region in which a solution of either Eq. \eqref{corrolary 1} or Eq. \eqref{corrolary 2} exists.
\end{corollary}
Let us now apply the mathematical tools presented in this section to specific examples.

\section{Examples \label{Examples section}}
\subsection{Kermack--McKendrick theory}
A simplified version of Kermack and McKendrick's mathematical theory of epidemics\cite{Kermack199133} is described by the following planar system with $x=(x_1,x_2)$\cite[see][p. 188]{strogatz2018nonlinear}:
\begin{eqnarray}
    \dot{x}=\begin{pmatrix}
        -k x_1 x_2 \\
    k x_1 x_2 -lx_2
    \end{pmatrix}. \label{Kermack system}
\end{eqnarray}
where $x_1\geq 0$ denotes the healthy population, $x_2\geq 0$ is the sick population and $k$, $l>0$ are constants. While we note that the Hamiltonian structure of a similar, three-dimensional model has previously been analyzed,\cite{Nutku1990L1145,Ballesteros2020} we focus here on the planar version. In the positive quadrant $x_1$, $x_2>0$, by \rr{comparison with Eq. \eqref{Transformed Hamiltonian system}}, Eq. \eqref{Kermack system} represents a transformed Hamiltonian system of the form of Eq. \eqref{Transformed Hamiltonian system} with 
\begin{eqnarray}
    \sqrt{\det g(x)}^{-1}&=&k x_1 x_2, \label{singular gdet}\\
    \widetilde{\mathcal{H}}(x)&=&\dfrac{l}{k}\ln x_1 -x_1-x_2. \label{KK Hamiltonian}
\end{eqnarray}
By Eq. \eqref{singular gdet}, the determinant of the metric tensor is singular on the two axes. As a consequence, the system given by Eq. \eqref{Kermack system} can not be expressed in terms of the same transformed Hamiltonian $\widetilde{\mathcal{H}}$ throughout $\mathbb{R}^2$. 

\subsection{Coupled Kuramoto oscillators \label{Swarm system}}
\begin{figure*}[t!]
\begin{psfrags}
    
\psfrag{a}{}
\psfrag{b}{\hspace{1cm} $\mathcal{K}$}
\psfrag{c}{\lc{\hspace{-0cm}$-1$\hspace{2.85cm}$0$\hspace{3cm}$1$}}
\psfrag{d}{\lc{\hspace{-0.7cm}$-0.477$\hspace{2.75cm}$0$\hspace{3.15cm}$0.477$}}
\psfrag{e}{\hspace{1.3cm} $\log_{10} \mathcal{J}$ }
\psfrag{1}{\hspace{-0.09cm}$\theta$}
\psfrag{2}{\hspace{0.045cm}$y$}
\psfrag{3}{$t$}

    \centering
    \hspace{-0.05cm}\includegraphics[width=0.65\textwidth]{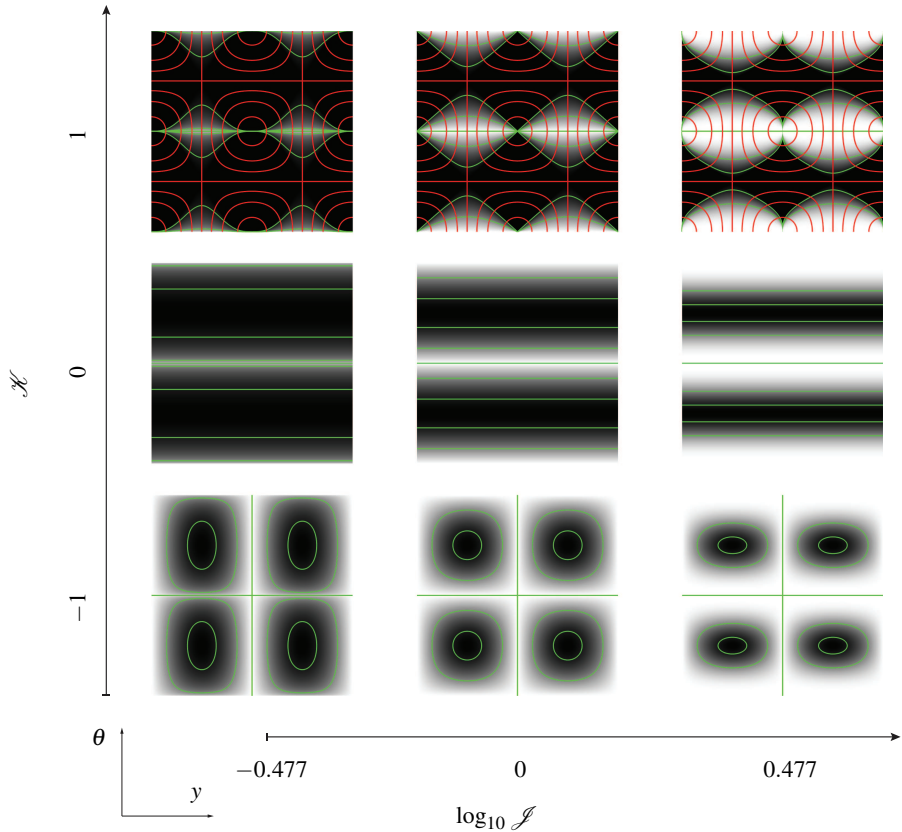}
\end{psfrags}
    \caption{Illustration of the Hamiltonian governing a system of coupled Kuramoto oscillators discussed in Sec. \ref{Swarm system}, over the coupling constants $\mathcal{J}$ and $\mathcal{K}$ \qq{(semi-log scale)}. Shown are contour plots of $-|\widetilde{\mathcal{H}}|$, where $\widetilde{\mathcal{H}}$ is defined by Eq. \eqref{swarm Hamiltonian}, over the periodic domain $[-\pi,\pi]\times [-\pi,\pi]$, ranging from values of $-1$ (black) to $0$ (white), overlaid with selected contour lines (green). In the cases with \qq{$\mathcal{K}=1$}, the colormap was cut off at $\pm 1$ because $\widetilde{\mathcal{H}}$ shows singular behavior. \mm{In the top row of this figure, the isocontours of the scalar potential $\widetilde{\mathcal{V}}=-\cos(y)\cos(\theta)$ derived in prior work,\cite{Pedergnana2022} governing the same system for $\mathcal{J},\mathcal{K}>0$, are overlaid in red over the isocontours of $\widetilde{\mathcal{H}}$.}}
    \label{Figure 3}
\end{figure*}

\citet{okeefe} derive, from two coupled Kuramoto oscillators, the following system, with $x=(y,\theta)$:
\begin{eqnarray}
  \dot{x}=\begin{pmatrix}
         -\mathcal{J}\sin y \cos  \theta,\\
     -\mathcal{K}\sin \theta \cos y, 
     \end{pmatrix} \label{swarm system new}
\end{eqnarray}
where $\mathcal{J}$ and $\mathcal{K}$ are coupling constants \mm{(see pp. 7--8 of the reference)}. We showed in our previous work\cite{Pedergnana2022} that this steady, planar system can be written as a transformed gradient system $\dot{x}=-g^{-1}\nabla \widetilde{\mathcal{V}}$ for $\mathcal{J}$, $\mathcal{K}>0$. If, however, $\mathcal{K}$ is less than zero, the same system has nested, closed orbits and there exists a conserved Hamiltonian $\widetilde{\mathcal{H}}(x)$.\cite{okeefe} In the positive quadrant $y$, $\theta>0$, Eq. \eqref{swarm system new} is equivalent to a transformed Hamiltonian system of the form of Eq. \eqref{Transformed Hamiltonian system} with\cite[see][see p. 7 for the derivation of $\widetilde{\mathcal{H}}$]{okeefe}
\begin{eqnarray}
\sqrt{\det g(x)}^{-1}&=&\dfrac{(\sin y)^{\mathcal{K}+1}}{(\sin \theta)^{\mathcal{J}-1}}, \label{swarm determinant}\\
    \widetilde{\mathcal{H}}(x)&=&-\dfrac{(\sin \theta)^\mathcal{J}}{(\sin y)^\mathcal{K}}. \label{swarm Hamiltonian}
\end{eqnarray}
Similar expressions can be found in each of the other quadrants \qq{from the requirement that $\det g$ is always positive}. \qq{A nonsmooth Hamiltonian which describes the dynamics governed by Eq. \eqref{swarm system new} in each of the quadrants (but not along the axes) is given by $-|\widetilde{\mathcal{H}}|$, where $\widetilde{\mathcal{H}}$ is defined in Eq. \eqref{swarm Hamiltonian}. This quantity is visualized in Fig. \ref{Figure 3} for different values of $\mathcal{J}$ and $\mathcal{K}$. Note that the color map is cut off in the top row of Fig. \ref{Figure 3}. This is because the main emphasis of this figure is the \textit{shape} of the isocontours, not their level values.}

The Hamiltonian formulation of Eq. \eqref{swarm system new} \qq{given in Eqs. \eqref{swarm determinant} and \eqref{swarm Hamiltonian}} \qq{is also true} for $\mathcal{K}>0$, where the system is known to be governed by \qq{the gradient of} a steady potential $\widetilde{\mathcal{V}}=-\cos(y)\cos(\theta)$ \mm{and the metric tensor $g=\mathrm{diag}(\mathcal{J},\mathcal{K})$.\cite{Pedergnana2022}} \qq{In other words, for $\mathcal{K}>0$, the system defined by Eq. \eqref{swarm system new} is simultaneously a Hamiltonian system and a gradient system.} Consequently, in this parameter range, $\widetilde{\mathcal{H}}$ can have no closed isocontours, because this property would be in contradiction to its gradient-driven nature for positive $\mathcal{J}$ and $\mathcal{K}$.\cite[see][pp. 201--202]{strogatz2018nonlinear} This argument is confirmed by the visualization of the Hamiltonian $\widetilde{\mathcal{H}}$ in Fig. \qq{\ref{Figure 2}}. This figure shows that, \qq{as expected, when} $\mathcal{K}$ is varied from negative to positive values, the contours of $\widetilde{\mathcal{H}}$ lose their closed character. \mm{In the top row of Fig. \ref{Figure 3}, the isocontours of the scalar potential $\widetilde{\mathcal{V}}$ are overlaid in red over the isocontours of $\widetilde{\mathcal{H}}$. Similar to the example in Sec. \eqref{prior work section}, for $\mathcal{J}=\mathcal{K}=1$ (top row, middle inset of Fig. \ref{Figure 3}), $g$ equals the identity matrix and the isocontours of $\widetilde{\mathcal{V}}$ are perpendicular to those of $\widetilde{\mathcal{H}}$.}

\subsection{Strogatz's conservative system \label{Strogatz system}}
Strogatz proposes the following conservative system with $x=(x_1,x_2)$:\cite[see][p. 190]{strogatz2018nonlinear}
\begin{eqnarray}
    \dot{x}=\begin{pmatrix}
        x_1 x_2 \\
        -x_1^2
    \end{pmatrix}. \label{Strogatz system new}
\end{eqnarray}
By \mm{comparison with Eq. \eqref{Transformed Hamiltonian system}}, in the right half-plane $x_1>0$, this system is a transformed Hamiltonian system \rr{with}
\begin{eqnarray}
\sqrt{\det g(x)}^{-1}&=&x_1,\\
    \widetilde{\mathcal{H}}(x)&=&\dfrac{x_1^2+x_2^2}{2}.
\end{eqnarray}
A similar expression can be found for the left half-plane. 

\subsection{Limit cycle scattering \label{SLS section}}
\begin{figure}[t]
\begin{psfrags}
\psfrag{1}{$A$}
\psfrag{2}{$\Phi$}
\psfrag{3}{$t$}
\psfrag{d}{\hspace{0.25cm}$\Delta=0$\hspace{1cm}$\Delta=40\pi$  rad/s\hspace{1cm}$\Delta=80\pi$ rad/s}

\psfrag{c}{\textcolor{red}{$\widetilde{\mathcal{V}}$}}

\begin{center}
\includegraphics[width=0.45\textwidth]{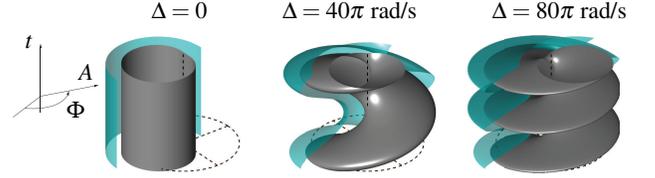}
\end{center}
\end{psfrags}
\vspace{-0.5cm}
    \caption{\mm{Isosurfaces of the unsteady potential \eqref{V eq 2 2} governing superradiant scattering by a limit cycle.\cite{pedergnana2023superradiant} The parameter values used to generate this figure correspond to those used in Fig. 3(a) in the reference. The variable parameter $\Delta=\omega_0-\omega$ denotes the detuning between the incident wave and the eigenfrequency of the limit cycle. The radius of the dashed circle corresponds to $10 A_0$, where $A_0$ is the unforced limit cycle amplitude. The dashed vertical line's length is $(\beta-\gamma)^{-1}$. The gray and cyan surfaces correspond to $0.5\%$ and $5\%$ of the maximum value of $\widetilde{\mathcal{V}}$ over the domain $[-20 A_0,20 A_0]\times[-20 A_0,20 A_0]\times [0,(\beta-\gamma)^{-1}]$. The explicit time-dependence visible in this figure arises from the interplay between the Hamiltonian and gradient parts of the Helmholtz decomposition in a steady system. }}
    \label{Figure 4}
\end{figure}

\mm{We consider a stochastic extension of the model introduced in \citet{pedergnana2023superradiant} to model \mm{superradiant scattering by a limit cycle}:
\begin{eqnarray}
\dot{a}=\Big(i\omega_0 -\gamma a+\dfrac{\beta}{1+\kappa |a|^2}\Big)+D s e^{i \omega t}+Z, \label{modal dynamics example}
\end{eqnarray}
where $a$ is the complex modal amplitude, $i$ is the imaginary unit, $\omega_0$ is the eigenfrequency, $\gamma$ is the damping, $\beta$ is the linear gain, $\kappa$ is the saturation constant, $D_j$ is the $j$\textsuperscript{th} entry of the coupling matrix $D$, $s$ is the incident wave amplitude, $\omega$ is the frequency of the incident wave, $Z=\xi_1+i\xi_2$, and $\xi_{1,2}$ is are white Gaussian noise sources. Defining $A=|a|$, $\phi=\arg (a)$, $\varphi=\phi-\omega t$ and separating real and imaginary parts yields, with $x=(A,\varphi)^T$,
\begin{eqnarray}
\dot{A}&=&-\gamma A+\frac{\beta A}{1+\kappa A^2}+\vert D_{j} \vert s \cos{(\arg D_{j}+\varphi)}+\xi_1, \nonumber \\
\label{SLS 1}\\
    \dot{\varphi}&=&\Delta-\frac{\vert D_{j} \vert s\sin{(\arg D_{j}+\varphi)}}{ A}+\frac{\xi_2}{A},\label{SLS 2}
\end{eqnarray}
where $\Delta=\omega_0-\omega$ is the detuning. In the absence of an incident wave ($s=0$), Eqs. \eqref{SLS 1} and \eqref{SLS 2} describe a stable limit cycle $a_0(t)=A_0 e^{i\omega_0 t}$ with constant amplitude $A_0=\sqrt{(\beta/\gamma-1)/\kappa}$. Assuming a curvilinear basis with $g=\mathrm{diag}(1,A^2)$ and $\sqrt{\det g}=A$ (see also Appendix \ref{Example Harmonic Oscillator}), a transformed Helmholtz decomposition \rr{in the form of Eq.} \eqref{Transformed 2D Helmholtz decomposition} of the system given by Eq. \eqref{modal dynamics example} is readily found to be 
\begin{eqnarray}
   \widetilde{\mathcal{V}}(x)&=&\dfrac{\gamma A^2}{2}-\dfrac{\beta \log(\kappa A^2+1)}{2 \kappa}\nonumber\\
   &&-|D_j| s A \cos(\varphi+\arg(D_j)),\label{V eq 2}\\
    \widetilde{\mathcal{H}}(x)&=&\dfrac{\Delta A^2}{2}, \label{H eq 2}
\end{eqnarray}
with the \rr{inverse} metric tensor \rr{$g^{-1}=\mathrm{diag}(1,A^{-2})$} and the positive definite matrix $h^{-1}=\mathrm{diag}(1,A^{-1})$ multiplying the noise vector $\Xi$. Unless $\omega$ is set to $\omega_0$ (no detuning), which makes the system given by Eq. \eqref{modal dynamics example} purely driven by the steady scalar potential, the interplay between the Hamiltonian and gradient parts of the Helmholtz decomposition is not easy to grasp. In contrast, the following alternative formulation with pure, but time-dependent, gradient is straightforward to interpret. This alternative description is obtained by absorbing the detuning $\Delta$ into the shifted phase $\Phi=\varphi-\Delta t$, leading to:
\begin{eqnarray}
   \widetilde{\mathcal{V}}(x,t)&=&\dfrac{\gamma A^2}{2}-\dfrac{\beta \log(\kappa A^2+1)}{2 \kappa}\nonumber\\
   &&-|D_j| s A \cos(\Phi+\Delta t+\arg(D_j)),\label{V eq 2 2}\\
    \widetilde{\mathcal{H}}(x)&=&0, \label{H eq 2 2}
\end{eqnarray}
with $x=(A,\Phi)^T$ and $g=h^T h=\mathrm{diag}(1,A^{-2})$. Equations \eqref{V eq 2 2} and \eqref{H eq 2 2} also describe a Helmholtz decomposition \eqref{Transformed 2D Helmholtz decomposition} of the system defined by Eq. \eqref{modal dynamics example}. As discussed in prior work on a similar example,\cite[see][Sec. V. A.]{Pedergnana2022} the explicit and periodic time-dependence of the potential $\widetilde{\mathcal{V}}$ in Eq. \eqref{V eq 2 2} corresponds to beating oscillations of $A$ and $\varphi$ (see Fig. \eqref{Figure 4}. For small detuning $\Delta\ll\omega$, these oscillations are slow compared to the forcing by the incident wave, and the system can be approximated as perfectly synchronized to good accuracy. The unsteady potential formulation given in Eqs. \eqref{V eq 2 2} and \eqref{H eq 2 2} reveals explicitly the time-dependent, deterministic forcing of the amplitude-phase dynamics for nonzero detuning, which is not directly evident from the steady, mixed Helmholtz decomposition given by Eqs. \eqref{V eq 2} and \eqref{H eq 2}. Potential future applications of the unsteady potential formulation are discussed in Sec. \ref{Discussion section}.}

\mm{We note that no term proportional to $A^{-1}$ appears in the modal dynamics defined by Eq. \eqref{modal dynamics example}, which would be the case had these equations been consistently derived by deterministic\cite{sanders2007averaging} and stochastic averaging\cite{stratonovich1963topics,Roberts1986111} of a corresponding ``fast oscillating'' system (see, for instance, \citet{Noiray2017}). Such a fast system may involve more intricate synchronization dynamics.\cite{balanov2009simple} Exploring the equivalence class of systems leading to Eq. \eqref{modal dynamics example} via the averaging method, also in view of their Helmholtz decomposition \eqref{Transformed 2D Helmholtz decomposition}, is a topic for future research.}

\subsection{Harmonic oscillator \label{Ruling out HO section}}
\begin{figure}[t!]
\begin{psfrags}
    
\psfrag{a}{\hspace{-0.1cm}$-4$}
\psfrag{c}{$0$}
\psfrag{e}{\hspace{-0.1cm}$4$}
\psfrag{k}{\hspace{-0.1cm}$-2$}
\psfrag{m}{$0$}
\psfrag{o}{\hspace{-0.1cm}$2$}
\psfrag{j}{\hspace{-0.1cm}$4$}
\psfrag{i}{}
\psfrag{h}{\hspace{-0.07cm}$0$}
\psfrag{g}{}
\psfrag{f}{\hspace{-0.37cm}$-4$}
\psfrag{t}{\hspace{-0.1cm}$2$}
\psfrag{s}{}
\psfrag{r}{\hspace{-0.12cm}$0$}
\psfrag{q}{}
\psfrag{p}{\hspace{-0.35cm}$-2$}
\psfrag{u}{\hspace{-0.5cm}}
\psfrag{v}{\hspace{-0.7cm}}
\psfrag{w}{\hspace{-0.55cm}}

\psfrag{A}{$x_1$}
\psfrag{B}{$x_2$}

    \centering
    \hspace{-0.05cm}\includegraphics[width=0.45\textwidth]{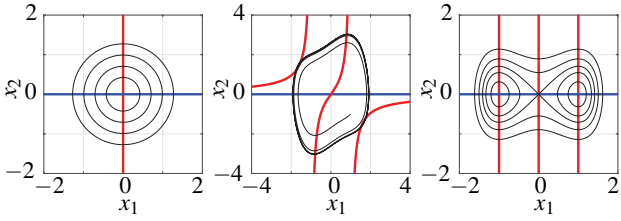}
\end{psfrags}
    \caption{Illustration of the results derived in \mm{Sec. \ref{Ruling out HO section}, Sec. \ref{Vdp section}, and Sec. \ref{Duffing section} for the Harmonic (left), Van der Pol (middle) and Duffing oscillators (right), respectively.}  By Theorem \ref{Theorem 1}, no closed orbits can be fully contained in any simply connected region throughout which a positive definite matrix function $\mathcal{N}$ exists \qq{such that Eq. \eqref{criterion theorem 1} is identically satisfied in that domain.} Using Corollary \ref{Corollary 1}, in Cartesian coordinates, two matrix functions satisfying this condition can be found by solving the first-order differential equations \eqref{corrolary 1} and \eqref{corrolary 2} for their respective off-diagonal elements. These solutions do \textit{not} exist, but are singular, over the blue and red curves, respectively, shown in the \mm{figure}. Therefore, any closed orbits of \qq{the systems \eqref{Unforced Damped harmonic oscillator},} \eqref{VDP system} and \eqref{Duffing system} must cross at least one blue and one red curve. Typical trajectories of \qq{all three} systems are shown in black, demonstrating that Theorem \ref{Theorem 1} is indeed satisfied for these examples. In the middle inset, $\alpha=0.7$.}
    \label{Figure 5}
\end{figure}

To apply Theorem \ref{Theorem 1} and Corollary \ref{Corollary 1}, we now consider the harmonic oscillator \qq{with $\omega_0=1$} in the \qq{undamped ($\gamma=0$),} unforced ($F=0$), noise-free ($\Xi=0$) limit, which is a steady planar \qq{Hamiltonian system $\dot{x}=S\nabla \mathcal{H}$ with $\mathcal{H}(x)=(x_1^2+x_2^2)/2$ that possesses nested closed orbits in the form of circles \qq{(or ellipses, if the coordinates are not normalized)} around the origin (see Fig. \ref{Figure 3}, left inset)}. In a Cartesian basis with \qq{normalized variables} $x=(x_1,x_2)$ \qq{(see Sec. \ref{Example Harmonic Oscillator})}, the \qq{dynamics} read 
\begin{eqnarray}
    \dot{x}&=&\underbrace{\begin{pmatrix}  x_2 \\
    -x_1
    \end{pmatrix}}_{=u(x)}. \label{Unforced Damped harmonic oscillator}
    \end{eqnarray}
    Substituting the velocity field $u(x)$ given by the RHS of Eq. \eqref{Unforced Damped harmonic oscillator} into Eqs. \eqref{corrolary 1} and \eqref{corrolary 2} yields the following set of differential equations:
        \begin{eqnarray}
        -a(x_1)+\dfrac{\partial \mathcal{N}_{12}(x_1,x_2)}{\partial x_2} x_1-b  =0, \label{DHO corrolary 1}\\
        -c+\dfrac{\partial \mathcal{N}_{21}(x_1,x_2)}{\partial x_1} x_2 - d(x_2)=0. \label{DHO corrolary 2}
    \end{eqnarray}
The general solutions of Eqs. \eqref{DHO corrolary 1} and \eqref{DHO corrolary 2} are given by
\begin{eqnarray}
     \mathcal{N}_{12}(x_1,x_2)&=&\dfrac{[a(x_1)+b] x_2}{x_1}+C_1, \label{DHO sol 1}\\
     \mathcal{N}_{21}(x_1,x_2)&=&\dfrac{[c+d(x_2)] x_1}{x_2}+C_2, \label{DHO sol 2}
\end{eqnarray}
where $C_1$ and $C_2$ are constants of integration. The \qq{expression} \eqref{DHO sol 1} exists for \qq{$x_1\neq 0$ (left and right half-plane)} and \eqref{DHO sol 2} for \qq{$x_2\neq 0$ (bottom and upper half-plane)}. None of \qq{these} regions fully contain any closed orbits of the system defined by Eq. \eqref{Unforced Damped harmonic oscillator}, confirming Theorem \ref{Theorem 1} and its Corollary \ref{Corollary 1}. \qq{The above results are visualized, together with typical trajectories of Eq. \eqref{Unforced Damped harmonic oscillator}, in Fig. \rr{\ref{Figure 5}} (left inset).}

\subsection{Van der Pol oscillator \label{Vdp section}}
Nearly a century ago, Balthasar Van der Pol proposed his nonlinear oscillator to describe self-sustained relaxation oscillations in electronic circuits.\cite{van1926lxxxviii} To test Theorem \ref{Theorem 1}, we consider here the following nondimensionalized set of equations describing an unforced Van der Pol oscillator\cite[see][pp. 67--68]{guckenheimer2013nonlinear}:
\begin{eqnarray}
    \dot{x}=\begin{pmatrix}
        x_2\\
        -x_1 +\alpha (1-x_1^2) x_2, 
    \end{pmatrix} \label{VDP system}
\end{eqnarray}
where $\alpha\geq 0$ is a non-negative constant. Assuming a Cartesian basis and substituting the velocity field $u(x)$ given by the RHS of \eqref{VDP system} into Eqs. \eqref{corrolary 1} and \eqref{corrolary 2} yields the following equations:
    \begin{eqnarray}
        -a(x_1)-\dfrac{\partial \mathcal{N}_{12}(x_1,x_2)}{\partial x_2}(-x_1 +\alpha (1-x_1^2) x_2)\nonumber\\
        -\mathcal{N}_{12}(x_1,x_2)\alpha(1-x_1^2)-b(1+2  \alpha x_1 x_2)=0, \label{VDP corrolary 1}\\
        -c+\dfrac{\partial \mathcal{N}_{21}(x_1,x_2)}{\partial x_1}x_2-d(x_2)(1+2\alpha x_1 x_2)=0. \label{VDP corrolary 2}
    \end{eqnarray}
These equations have the following general solutions:
\begin{eqnarray}
    \mathcal{N}_{12}(x_1,x_2)&=&\dfrac{[a(x_1)+b]x_2+ \alpha  b x_1 x_2^2+C_1}{x_1 -\alpha (1-x_1^2) x_2}, \label{VDP sol 1}\\
    \mathcal{N}_{21}(x_1,x_2)&=&\dfrac{[c+d(x_2)]x_1}{x_2}+ d(x_2) \alpha x_1^2 +C_2. \label{VDP sol 2}
\end{eqnarray}
For $\alpha>0$, the zero level set of the denominator of \eqref{VDP sol 1} describes a curve which crosses the origin and diverges to negative infinity at $x_1=-1$ and to positive infinity at $x_1=1$. In general, $\mathcal{N}_{12}$ exists everywhere but on this curve. Similarly, the solution \eqref{VDP sol 2} exists in the upper and in the lower half-plane, but not on the $x$-axis ($x_2=0$). The limit cycle of the Van der Pol oscillator,\cite[see][pp. 67--82]{guckenheimer2013nonlinear} which winds around the origin, necessarily crosses both of these curves, on which the solutions \eqref{VDP sol 1} and \eqref{VDP sol 2} do not exist, confirming Theorem \ref{Theorem 1}. See Fig. \rr{\ref{Figure 5}}, middle inset, for an illustration of the above results on an example with $\alpha=0.7$.

\subsection{Duffing nonlinearity \label{Duffing section}}
For our next example, we examine the Duffing oscillator\cite{duffing1918forced} with a double-well Hamiltonian in nondimensionalized form\cite[see][pp. 82--91]{guckenheimer2013nonlinear}:
\begin{eqnarray}
    \dot{x}=\begin{pmatrix}
        x_2\\
        x_1-x_1^3. 
    \end{pmatrix} \label{Duffing system}
\end{eqnarray}
Equation \eqref{Duffing system} is of the form $\dot{x}=S\nabla \mathcal{H}$ with $\mathcal{H}(x)=(x_2^2-x_1^2)/2+\rr{x_1}^4/4$ and is governed by two families of closed orbits, segregated by the curve $\mathcal{H}(x)=0$, the separatrix. For $\mathcal{H}>0$, the flow \eqref{Duffing system} has closed orbits encircling the separatrix and the fixed point at the origin. In the region where $\mathcal{H}<0$, inside the separatrix, there are two nested sets of closed orbits, mirror-symmetric with respect to the $y$-axis, which wind around the fixed points $x=(-1,0)^T$ and $x=(1,0)^T$, respectively (see Fig. \rr{\ref{Figure 5}}, right inset). Using a Cartesian basis, substituting the RHS of \eqref{Duffing system} into Eqs. \eqref{corrolary 1} and \eqref{corrolary 2} leads to
       \begin{eqnarray}
        -a(x_1)-\dfrac{\partial \mathcal{N}_{12}(x_1,x_2)}{\partial x_2}(x_1-x_1^3)+b(1-3x_1^2)=0, \label{Duffing corrolary 1}\\
        -c+\dfrac{\partial \mathcal{N}_{21}(x_1,x_2)}{\partial x_1}x_2+d(x_2)(1-3x_1^2)=0, \label{Duffing corrolary 2}
    \end{eqnarray}
whose general solutions are given by
\begin{eqnarray}
    \mathcal{N}_{12}(x_1,x_2)&=&\dfrac{[b-a(x_1)]x_2-3bx_1^2 x_2}{x_1-x_1^3}+C_1, \label{Duffing sol 1}\\
    \mathcal{N}_{21}(x_1,x_2)&=&\dfrac{c x_1+d(x_2)(x_1^3-x_1)}{x_2}+C_2. \label{Duffing sol 2}
\end{eqnarray}
By Corollary \ref{Corollary 1}, looking at Eqs. \eqref{Duffing sol 1} and \eqref{Duffing sol 2}, any closed orbits of the steady planar system defined by Eq. \eqref{Duffing system} must cross one of the lines \rr{$x_1=1$}, $-1$ or \rr{$x_1=0$} and, additionally, cross the $x$-axis $x_2=0$. By the above discussion, as shown in Fig. \rr{\ref{Figure 5}}, right inset, these conditions are indeed satisfied, confirming Theorem \ref{Theorem 1}.

\subsection{Quadratic system with four limit cycles \label{Shi section}}
\begin{figure}[t!]
\begin{psfrags}
    
\psfrag{a}{\hspace{-0.2cm}$-2$\hspace{0.8cm}$-1$}
\psfrag{c}{\hspace{-0.05cm}$0$\hspace{1.03cm}$1$}
\psfrag{e}{\hspace{-0.02cm}$2$}
\psfrag{f}{\hspace{-0.45cm}$-2$}
\psfrag{g}{\hspace{-0.43cm}$-1$}
\psfrag{h}{\hspace{-0.15cm}$0$}
\psfrag{i}{\hspace{-0.2cm}$1$}
\psfrag{j}{\hspace{-0.2cm}$2$}
\psfrag{x}{\hspace{-0.05cm}$x_1$}
\psfrag{k}{$x_2$}

  \includegraphics[width=0.3\textwidth]{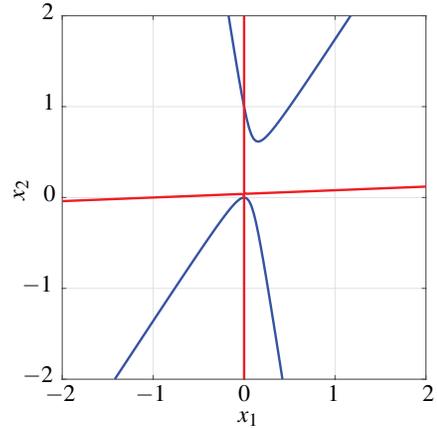}
\end{psfrags}
    \caption{\mm{Singularity curves of the solutions $\mathcal{N}_{12}$ and $\mathcal{N}_{21}$ of Eqs. \eqref{corrolary 1} (red) and \eqref{corrolary 2} (blue), respectively, for the system of Shi Songling.\cite{shi1980concrete,Songling1981301} Three of the four limit cycles of this system (not shown here, see Fig. 2 of \citet{Galias2022}) are tiny and encircle the point $(0,0)$. The fourth limit cycle encircles the point $(0,1)$. For a given system, \rr{similar curves as shown in this figure} could be used as one-dimensional manifolds to seed initial conditions from for brute-force detection of limit cycles. In such an algorithm, the distances between intersections and critical points of the curves could serve as natural length scales for determining the discretization fineness in different regions.}}
    \label{Figure 6}
\end{figure}

\mm{We consider the system of Shi Songling\cite{shi1980concrete,Songling1981301}, given by\cite{Galias2022}
\begin{eqnarray}
    \dot{x}=\begin{pmatrix}
        \lambda x_1-x_2-10 x_1^2+(5+\delta) x_1 x_2+x_2^2\\
x_1+x_1^2+(-25+8\varepsilon-9\delta) x_1 x_2
    \end{pmatrix}, \label{Shi system}
\end{eqnarray}
for the parameter values $\delta=-10^{-13}$, $\varepsilon=-10^{-52}$ and 
$\lambda=-10^{-200}$. This system has exactly four limit cycles. Three of the limit cycles are tiny (see Fig. 2 of \citet{Galias2022}) and encircle the point $(0,0)$, while the fourth is normal-sized and encircles $(0,1)$. The interested reader is referred to \citet{Kuznetsov201329} for a quadratic planar system with four normal-sized limit cycles. The solutions of Eqs. \eqref{corrolary 1} and \eqref{corrolary 2} for the system given by Eq. \eqref{Shi system} are of the form
\begin{eqnarray}
    \mathcal{N}_{12}(x_1,x_2)&=&\frac{X(x_1,x_2)+C_1}{x_1+x_1^2+(-25+8\varepsilon-9\delta) x_1 x_2},
\\
        \mathcal{N}_{21}(x_1,x_2)&=&\dfrac{Y(x_1,x_2)+C_2}{ \lambda x_1-x_2-10 x_1^2+(5+\delta) x_1 x_2+x_2^2},
\end{eqnarray}
where\rr{, for suitable choices of $a(x_1)$ and $d(x_2)$,} $X$ and $Y$ are compositions of smooth, elementary functions. The curves on which $\mathcal{N}_{12}$ and $\mathcal{N}_{21}$ are singular are visualized in Fig. \ref{Figure 6}. These curves intersect at the points $(0,0)$ and $(0,1)$. We observe that the distance between the intersection at the origin and a nearby critical point is small. In a (semi-)automated analysis, this distance could be used as a first estimate of the length scale at which the dynamics near the origin take place, determining the discretization fineness of initial conditions seeded along the singularity curves of $\mathcal{N}_{12}$ and $\mathcal{N}_{21}$. }

\section{Discussion and Outlook \label{Discussion section}}
\mm{This work explores smooth transformations of planar systems. It was shown in Sec. \ref{prior work section} and in Sec. \ref{Swarm system} that Hamiltonian systems without closed isocontours can simultaneously be gradient systems. In light of these results, it would be worth investigating if one could consistently define a gradient system from such a Hamiltonian and vice versa. If so, then this dualism may help towards quantization\cite{dirac1981principles} of systems which are not purely Hamiltonian, which is a topic of ongoing research with theoretical and practical significance.\cite{moiseyev2011non,Deguchi2020,Blacker2021,eichler2023classical} A geometric approach to unifying Hamiltonian and gradient dynamics is taken by \citet{Esen2022}. To analyze the saddle-type Hamiltonians discussed here, it may prove useful to introduce hyperbolic action-angle coordinates in the spirit of \citet{Waalkens2008R1} (see pp. 25--26).}

\mm{In Sec. \rr{\ref{SLS section}}, an example was given of a steady planar system with nonzero Hamiltonian and gradient part. It was shown that this system can alternatively be written as an \textit{unsteady}, purely gradient-driven system with vanishing Hamiltonian part, complementing earlier results.\cite{Pedergnana2022} It could be interesting to study whether such a reformulation is also possible for systems oscillating at more than one frequency.\cite{delpino2023limit}}

\mm{In Sec. \ref{Main results 2}, a criterion was presented for ruling out closed orbits in certain regions of phase space of steady systems. From this general result, a corollary was derived which can be applied in automated fashion, yielding ordinary differential equations which can be easily solved. In contrast, applying Dulac's criterion in automated fashion would yield (partial) differential inequalities, which are more complicated. As explained in Sec. \ref{Shi section}, these results could be used in the future for efficient and automated seeding of initial conditions in numerical algorithms to detect periodic solutions.\cite{parker2012practical,Jahn2020}}

\section*{Author contributions}
\cc{\textbf{Tiemo Pedergnana}: Conceptualization (lead), formal analysis, writing–original draft, writing–review and editing (lead). \textbf{Nicolas Noiray}: Conceptualization (supporting), supervision, writing–review and editing (supporting).}

\section*{Acknowledgements}
This project is funded by the Swiss National Science Foundation under Grant agreement 184617.

\section*{AUTHOR DECLARATIONS}
\subsection*{Conflict of Interest}
\dc{The authors have no conflicts to disclose.}
\section*{Data availability}
\dc{The datasets used for generating
the plots in this study can be directly obtained
by numerical simulation of the related mathematical equations in the manuscript.}

\appendix 

\section{Reformulation of Bendixson's criterion \label{Appendix A}}
\mm{Bendixson's criterion concerns the steady planar system $\dot{x}=u(x)$. The criterion states that, if the divergence of $u$ has a constant sign throughout a simply connected domain $\mathcal{D}$, then the system has no closed orbits fully contained in that domain. The coordinate-independent divergence operator for a vector field $u$ in curvilinear coordinates is given by \cite[see][p. 92]{jost2008riemannian}
\begin{eqnarray}
    \mathrm{div} \, u(x)=\dfrac{1}{\sqrt{\det g(x)}}\nabla \cdot \big[\sqrt{\det g(x)}\, u(x) \big].
\end{eqnarray}
Applying this operator to the RHS of Eq. \eqref{steady planar field} and requiring that the result be greater than zero yields 
\begin{eqnarray}
-\dfrac{1}{\sqrt{\det g(x)}}\nabla \cdot \big[\sqrt{\det g(x)}\,g^{-1}(x) \nabla \widetilde{\mathcal{V}}(x)\big]>0. \label{V inequality}
\end{eqnarray}
This requirement is simply the definition of a strictly subharmonic function \cite[see][p. 433]{jost2008riemannian}, and the operator acting on $\widetilde{\mathcal{V}}$ is the (negative) Laplace--Beltrami operator,\cite[see][p. 92]{jost2008riemannian} which arises here naturally from elementary manipulations. The analogous result for $-\widetilde{\mathcal{V}}$ is obtained by replacing ``$>$'' with ``$<$'' in the inequality \eqref{V inequality}.}

\section{Transformation formula for the Harmonic oscillator\label{Example Harmonic Oscillator} \label{Appendix B}}
\mm{In this appendix}, we consider the noise-driven, forced-damped harmonic oscillator with position $x_1$ and normalized velocity $x_2=\dot{x}_1/\omega_0$, whose dynamics are defined by the system
\begin{eqnarray}
    \dot{x}&=&\begin{pmatrix} \omega_0 x_2 \\
    -\gamma x_2-\omega_0x_1+\dfrac{F}{\omega_0}\cos\omega t+\dfrac{\xi}{\omega_0}
    \end{pmatrix} \label{Damped harmonic oscillator},
    \end{eqnarray}
where $\omega_0$ is the eigenfrequency, $\gamma$ is the damping, $F$ is the forcing amplitude and $\xi$ is a zero-mean Gaussian white noise source with variance $\Gamma$. The HD of the planar system defined by \rr{Eq.} \eqref{Damped harmonic oscillator} can be written as follows:
    \begin{eqnarray}
            \dot{x}&=&\underbrace{-\begin{pmatrix}
\frac{\partial}{\partial x_1}\\
\frac{\partial}{\partial x_2}
\end{pmatrix} \Big(\frac{\gamma x_2^2}{2}-\dfrac{F}{\omega_0}x_2\cos\omega t\Big)}_{=-\nabla \mathcal{V}(x,t)}+\underbrace{\begin{pmatrix}
\frac{\partial}{\partial x_2}\\
-\frac{\partial}{\partial x_1}
\end{pmatrix} \frac{\omega_0( x_1^2+x_2^2)}{2}}_{=S\nabla \mathcal{H}(x)} \nonumber \\
&&+\underbrace{\begin{pmatrix}
    0\\
    \dfrac{\xi}{\omega_0}
\end{pmatrix}}_{=\Xi}. \label{HD int res}
\end{eqnarray}
Note that the HD in Eq. \eqref{HD int res} is only quasi-unique, as the forcing term could be neglected in the potential $\mathcal{V}$ and, instead, a term proportional to $x_1$ (specifically, $-[F x_1/\omega_0] \cos\omega t$) could be added to the Hamiltonian $\mathcal{H}$, leaving Eq. \eqref{Damped harmonic oscillator} unchanged.

Next, we apply a transformation to amplitude-phase coordinates $y=(A,\phi)$, $A\in\mathbb{R}^+$, $\phi\in [0, 2\pi]$:
\begin{eqnarray}
    \underbrace{\begin{pmatrix}
x_1\\
x_2
\end{pmatrix}}_{x}= \underbrace{\begin{pmatrix}
A \cos {\phi}\\
 A\sin {\phi}
\end{pmatrix}}_{f(y)}. \label{DHO transformation}
\end{eqnarray}
After substituting the expression \eqref{DHO transformation} into Eq. \eqref{Damped harmonic oscillator}, some algebraic manipulation and redefining $y\rightarrow x$, we find that the transformed dynamics are given by, with $x=(A,\phi)$:
\begin{eqnarray}
   \dot{x}&=& \begin{pmatrix}
-\gamma A \sin^2\phi+\dfrac{F}{\omega_0}\cos\omega t\sin\phi+\dfrac{\xi}{\omega_0}\sin \phi\\
-\omega_0- \gamma \sin\phi\cos\phi+\dfrac{F}{\omega_0 A}\cos\omega t \cos\phi +\dfrac{\xi}{\omega_0 A}\cos\phi
\end{pmatrix}\nonumber \\\label{DHO tf dynamics}
\end{eqnarray}
To confirm the transformation formula in Eq. \eqref{Transformed 2D Helmholtz decomposition}, we compute the Jacobian matrix of the mapping \eqref{DHO transformation} and its polar decomposition:
\begin{eqnarray}
    J(x)&=&\underbrace{\begin{pmatrix}
\cos \phi & -\sin\phi \\
 \sin\phi &  \cos \phi
\end{pmatrix}}_{=Q(x)}\underbrace{\begin{pmatrix}
1 & 0 \\
0 & A
\end{pmatrix}}_{=h(x)}. \label{J Ex A}
\end{eqnarray}
We observe that $\det Q=1$. From the expressions in Eq. \eqref{J Ex A}, we immediately obtain the metric tensor $g=h^T h$ and its determinant\qq{:}
\begin{eqnarray}
g(x)&=&\mathrm{diag}(1,A^2),\\
\det g(x)&=&A^2.
\end{eqnarray}
Furthermore, \mm{we have} $\widetilde{\Xi}=Q^T\Xi=(\xi \sin\phi/\omega_0,\xi\cos\phi/\omega_0)^T$. By substituting the above expressions into Eq. \eqref{Transformed 2D Helmholtz decomposition}, using the definitions $\widetilde{\mathcal{V}}(x,t)=\mathcal{V}(f(x),t)$ and $\widetilde{\mathcal{H}}(x,t)=\mathcal{H}(f(x),t)$, one can verify that the transformed dynamics, defined by \rr{Eq.} \eqref{DHO tf dynamics}, are equivalent to
\begin{eqnarray}
    \dot{x}&=&-\underbrace{\begin{pmatrix}
1 & 0 \\
0 & A^{-2}
\end{pmatrix}}_{=g^{-1}(x)}\underbrace{\begin{pmatrix}
\frac{\partial}{\partial A}\\
\frac{\partial}{\partial \phi}
\end{pmatrix}\bigg(\dfrac{\gamma A^2 \sin^2\phi}{2}-\dfrac{F A}{\omega_0}\cos\omega t \sin\phi\bigg)}_{=\nabla\widetilde{\mathcal{V}}(x,t)}\nonumber\\
&+&\hspace{-0.1cm}\underbrace{\dfrac{1}{A}}_{=\sqrt{\det g(x)}^{-1}}\underbrace{\begin{pmatrix}
\frac{\partial}{\partial \phi}\\
-\frac{\partial}{\partial A}
\end{pmatrix}\bigg(\dfrac{\omega_0 A^2}{2}\bigg)}_{=S\nabla \widetilde{\mathcal{H}}(x)}+\underbrace{\begin{pmatrix}
1 & 0 \\
0 & A^{-1}
\end{pmatrix}}_{=h^{-1}(x)}\underbrace{\begin{pmatrix}
    \widetilde{\xi}_1\\
    \widetilde{\xi}_2
\end{pmatrix}}_{=\widetilde{\Xi}}, \label{DHO IHD}
\end{eqnarray}
confirming Eq. \eqref{Transformed 2D Helmholtz decomposition}.

\section{Polynomial Liénard systems \label{Appendix C}}
\qq{In first-order form, Liénard's equation reads as follows:\cite[see][p. 212]{strogatz2018nonlinear}
\begin{eqnarray}
    \dot{x}=\begin{pmatrix} x_2 \\
    -p(x_1)-q(x_1)x_2
    \end{pmatrix}. \label{Lienard system}
\end{eqnarray}
We focus here on the class of nonlinear oscillators for which the antiderivative of $p$ exists and $q$ is given by a polynomial of finite order $M$:
\begin{eqnarray}
    q(x_1)=\sum_{m=1}^M q_m x_1^m.\label{q eq} 
\end{eqnarray}
The subclass for which $M\leq 2$ encompasses many well-known examples of linear and nonlinear oscillators, \rr{some} of which are analyzed in \rr{this work}. The bifurcations of a system of the form of Eq. \eqref{Lienard system} with $M=2$ for which $p$ is a polynomial of order $5$ have recently been analyzed.\cite{delpino2023limit} A Liénard system, specified by Eq. \eqref{Lienard system} with $M=4$, has previously been studied in the context of thermoacoustic instabilities in turbulent combustors.\cite{Noiray2017}}

\qq{Assuming a Cartesian basis, the Helmholtz decomposition \eqref{Cartesian system} of the Liénard system governed by Eq. \eqref{Lienard system} is given by
\begin{eqnarray}
    \mathcal{V}(x)&=&\dfrac{q(x_1)x_2^2}{2} -q''(x_1)\dfrac{x_2^4}{24}+\dots\nonumber\\
    &=&\sum_{n=0} q^{(2n)}(x_1) \dfrac{(-1)^{n}x_2^{2(n+1)}}{2(n+1)!}, \label{V eq}\\
    \mathcal{H}(x)&=&\dfrac{x_2^2}{2}+\int p(x_1) dx_1+q'(x_1)\dfrac{x_2^3}{6}-q'''(x_1)\dfrac{x_2^5}{120}+\dots\nonumber\\
    &=&\dfrac{x_2^2}{2}+\int p(x_1) dx_1+\sum_{n=0} q^{(2n+1)}(x_1) \dfrac{(-1)^{n}x_2^{2n+3}}{(2n+3)!}, \nonumber\\ \label{H eq}
\end{eqnarray}
where a dash denotes the derivative, $\int p dx_1$ is the antiderivative of $p$ and $q^{(n)}$ denotes the $n$\textsuperscript{th} derivative of $q$. This result can be verified by substituting the above expressions into Eq. \eqref{Helmholtz decomposition}, whereby the Eq. \eqref{Lienard system} is recovered. The sums in Eqs. \eqref{V eq} and \eqref{H eq} truncate due to the finite order of $q$ defined in Eq. \eqref{q eq}.}

\bibliography{bibliography}

\begin{thebibliography}{84}%
\makeatletter
\providecommand \@ifxundefined [1]{%
 \@ifx{#1\undefined}
}%
\providecommand \@ifnum [1]{%
 \ifnum #1\expandafter \@firstoftwo
 \else \expandafter \@secondoftwo
 \fi
}%
\providecommand \@ifx [1]{%
 \ifx #1\expandafter \@firstoftwo
 \else \expandafter \@secondoftwo
 \fi
}%
\providecommand \natexlab [1]{#1}%
\providecommand \enquote  [1]{``#1''}%
\providecommand \bibnamefont  [1]{#1}%
\providecommand \bibfnamefont [1]{#1}%
\providecommand \citenamefont [1]{#1}%
\providecommand \href@noop [0]{\@secondoftwo}%
\providecommand \href [0]{\begingroup \@sanitize@url \@href}%
\providecommand \@href[1]{\@@startlink{#1}\@@href}%
\providecommand \@@href[1]{\endgroup#1\@@endlink}%
\providecommand \@sanitize@url [0]{\catcode `\\12\catcode `\$12\catcode `\&12\catcode `\#12\catcode `\^12\catcode `\_12\catcode `\%12\relax}%
\providecommand \@@startlink[1]{}%
\providecommand \@@endlink[0]{}%
\providecommand \url  [0]{\begingroup\@sanitize@url \@url }%
\providecommand \@url [1]{\endgroup\@href {#1}{\urlprefix }}%
\providecommand \urlprefix  [0]{URL }%
\providecommand \Eprint [0]{\href }%
\providecommand \doibase [0]{https://doi.org/}%
\providecommand \selectlanguage [0]{\@gobble}%
\providecommand \bibinfo  [0]{\@secondoftwo}%
\providecommand \bibfield  [0]{\@secondoftwo}%
\providecommand \translation [1]{[#1]}%
\providecommand \BibitemOpen [0]{}%
\providecommand \bibitemStop [0]{}%
\providecommand \bibitemNoStop [0]{.\EOS\space}%
\providecommand \EOS [0]{\spacefactor3000\relax}%
\providecommand \BibitemShut  [1]{\csname bibitem#1\endcsname}%
\let\auto@bib@innerbib\@empty
\bibitem [{\citenamefont {Stratonovich}(1963)}]{stratonovich1963topics}%
  \BibitemOpen
  \bibfield  {author} {\bibinfo {author} {\bibfnamefont {R.}~\bibnamefont {Stratonovich}},\ }\href@noop {} {\emph {\bibinfo {title} {Topics in the Theory of Random Noise Vol. I: General Theory of Random Processes Nonlinear Transformations of Signals and Noise}}}\ (\bibinfo  {publisher} {Gordon \& Breach},\ \bibinfo {year} {1963})\BibitemShut {NoStop}%
\bibitem [{\citenamefont {Gardiner}\ \emph {et~al.}(1985)\citenamefont {Gardiner} \emph {et~al.}}]{gardiner1985handbook}%
  \BibitemOpen
  \bibfield  {author} {\bibinfo {author} {\bibfnamefont {C.~W.}\ \bibnamefont {Gardiner}} \emph {et~al.},\ }\href@noop {} {\emph {\bibinfo {title} {Handbook of stochastic methods}}},\ Vol.~\bibinfo {volume} {4}\ (\bibinfo  {publisher} {Springer Berlin},\ \bibinfo {year} {1985})\BibitemShut {NoStop}%
\bibitem [{\citenamefont {Risken}(1984)}]{risken1996fokker}%
  \BibitemOpen
  \bibfield  {author} {\bibinfo {author} {\bibfnamefont {H.}~\bibnamefont {Risken}},\ }\href@noop {} {\emph {\bibinfo {title} {The Fokker-Planck Equation}}}\ (\bibinfo  {publisher} {Springer},\ \bibinfo {year} {1984})\BibitemShut {NoStop}%
\bibitem [{Note1()}]{Note1}%
  \BibitemOpen
  \bibinfo {note} {Throughout this work, a dot over a dependent variable denotes the total time derivative.}\BibitemShut {Stop}%
\bibitem [{\citenamefont {Guckenheimer}\ and\ \citenamefont {Holmes}(2013)}]{guckenheimer2013nonlinear}%
  \BibitemOpen
  \bibfield  {author} {\bibinfo {author} {\bibfnamefont {J.}~\bibnamefont {Guckenheimer}}\ and\ \bibinfo {author} {\bibfnamefont {P.}~\bibnamefont {Holmes}},\ }\href@noop {} {\emph {\bibinfo {title} {Nonlinear oscillations, dynamical systems, and bifurcations of vector fields}}},\ Vol.~\bibinfo {volume} {42}\ (\bibinfo  {publisher} {Springer {S}cience \& {B}usiness {M}edia},\ \bibinfo {year} {2013})\BibitemShut {NoStop}%
\bibitem [{\citenamefont {Strogatz}(2015)}]{strogatz2018nonlinear}%
  \BibitemOpen
  \bibfield  {author} {\bibinfo {author} {\bibfnamefont {S.~H.}\ \bibnamefont {Strogatz}},\ }\href@noop {} {\emph {\bibinfo {title} {Nonlinear dynamics and chaos: with applications to physics, biology, chemistry, and engineering}}}\ (\bibinfo  {publisher} {CRC press},\ \bibinfo {year} {2015})\BibitemShut {NoStop}%
\bibitem [{\citenamefont {Stokes}(2009)}]{stokes_2009}%
  \BibitemOpen
  \bibfield  {author} {\bibinfo {author} {\bibfnamefont {G.~G.}\ \bibnamefont {Stokes}},\ }\enquote {\bibinfo {title} {On the dynamical theory of diffraction},}\ in\ \href {https://doi.org/10.1017/CBO9780511702259.015} {\emph {\bibinfo {booktitle} {Mathematical and Physical Papers}}},\ \bibinfo {series} {Cambridge Library Collection - Mathematics}, Vol.~\bibinfo {volume} {2}\ (\bibinfo  {publisher} {Cambridge University Press},\ \bibinfo {year} {2009})\ p.\ \bibinfo {pages} {243–328}\BibitemShut {NoStop}%
\bibitem [{\citenamefont {Helmholtz}(1858)}]{Helmholtz185825}%
  \BibitemOpen
  \bibfield  {author} {\bibinfo {author} {\bibfnamefont {H.}~\bibnamefont {Helmholtz}},\ }\bibfield  {title} {\enquote {\bibinfo {title} {Über {I}ntegrale der {H}ydrodynamischen {G}leichungen, welche den {W}irbelbewegungen entsprechen},}\ }\href {https://doi.org/10.1515/crll.1858.55.25} {\bibfield  {journal} {\bibinfo  {journal} {J. Reine Angew. Math.}\ }\textbf {\bibinfo {volume} {1858}},\ \bibinfo {pages} {25--55} (\bibinfo {year} {1858})}\BibitemShut {NoStop}%
\bibitem [{\citenamefont {Morse}\ and\ \citenamefont {Feshbach}(1953)}]{morse1953methods}%
  \BibitemOpen
  \bibfield  {author} {\bibinfo {author} {\bibfnamefont {P.~M.}\ \bibnamefont {Morse}}\ and\ \bibinfo {author} {\bibfnamefont {H.}~\bibnamefont {Feshbach}},\ }\href@noop {} {\emph {\bibinfo {title} {Methods of Theoretical Physics. Vol. 1-2}}}\ (\bibinfo  {publisher} {New York},\ \bibinfo {year} {1953})\BibitemShut {NoStop}%
\bibitem [{\citenamefont {Arnold}(2013)}]{arnol2013mathematical}%
  \BibitemOpen
  \bibfield  {author} {\bibinfo {author} {\bibfnamefont {V.~I.}\ \bibnamefont {Arnold}},\ }\href@noop {} {\emph {\bibinfo {title} {Mathematical methods of classical mechanics}}},\ Vol.~\bibinfo {volume} {60}\ (\bibinfo  {publisher} {Springer Science \& Business Media},\ \bibinfo {year} {2013})\BibitemShut {NoStop}%
\bibitem [{\citenamefont {Morino}(1986)}]{Morino198665}%
  \BibitemOpen
  \bibfield  {author} {\bibinfo {author} {\bibfnamefont {L.}~\bibnamefont {Morino}},\ }\bibfield  {title} {\enquote {\bibinfo {title} {Helmholtz decomposition revisited: Vorticity generation and trailing edge condition - part 1: Incompressible flows},}\ }\href {https://doi.org/10.1007/BF00298638} {\bibfield  {journal} {\bibinfo  {journal} {Comput. Mech.}\ }\textbf {\bibinfo {volume} {1}},\ \bibinfo {pages} {65--90} (\bibinfo {year} {1986})}\BibitemShut {NoStop}%
\bibitem [{\citenamefont {Joseph}(2006)}]{Joseph200614272}%
  \BibitemOpen
  \bibfield  {author} {\bibinfo {author} {\bibfnamefont {D.}~\bibnamefont {Joseph}},\ }\bibfield  {title} {\enquote {\bibinfo {title} {Helmholtz decomposition coupling rotational to irrotational flow of a viscous fluid},}\ }\href {https://doi.org/10.1073/pnas.0605792103} {\bibfield  {journal} {\bibinfo  {journal} {Proc. Natl. Acad. Sci. U.S.A.}\ }\textbf {\bibinfo {volume} {103}},\ \bibinfo {pages} {14272--14277} (\bibinfo {year} {2006})}\BibitemShut {NoStop}%
\bibitem [{\citenamefont {Linke}(2014)}]{Linke2014782}%
  \BibitemOpen
  \bibfield  {author} {\bibinfo {author} {\bibfnamefont {A.}~\bibnamefont {Linke}},\ }\bibfield  {title} {\enquote {\bibinfo {title} {On the role of the {H}elmholtz decomposition in mixed methods for incompressible flows and a new variational crime},}\ }\href {https://doi.org/10.1016/j.cma.2013.10.011} {\bibfield  {journal} {\bibinfo  {journal} {Comput. Methods Appl. Mech. Eng.}\ }\textbf {\bibinfo {volume} {268}},\ \bibinfo {pages} {782--800} (\bibinfo {year} {2014})}\BibitemShut {NoStop}%
\bibitem [{\citenamefont {Bühler}, \citenamefont {Callies},\ and\ \citenamefont {Ferrari}(2014)}]{Buhler20141007}%
  \BibitemOpen
  \bibfield  {author} {\bibinfo {author} {\bibfnamefont {O.}~\bibnamefont {Bühler}}, \bibinfo {author} {\bibfnamefont {J.}~\bibnamefont {Callies}},\ and\ \bibinfo {author} {\bibfnamefont {R.}~\bibnamefont {Ferrari}},\ }\bibfield  {title} {\enquote {\bibinfo {title} {Wave-vortex decomposition of one-dimensional ship-track data},}\ }\href {https://doi.org/10.1017/jfm.2014.488} {\bibfield  {journal} {\bibinfo  {journal} {J. Fluid Mech.}\ }\textbf {\bibinfo {volume} {756}},\ \bibinfo {pages} {1007--1026} (\bibinfo {year} {2014})}\BibitemShut {NoStop}%
\bibitem [{\citenamefont {Lindborg}(2014)}]{Lindborg2014}%
  \BibitemOpen
  \bibfield  {author} {\bibinfo {author} {\bibfnamefont {E.}~\bibnamefont {Lindborg}},\ }\bibfield  {title} {\enquote {\bibinfo {title} {A {H}elmholtz decomposition of structure functions and spectra calculated from aircraft data},}\ }\href {https://doi.org/10.1017/jfm.2014.685} {\bibfield  {journal} {\bibinfo  {journal} {J. Fluid Mech.}\ }\textbf {\bibinfo {volume} {762}} (\bibinfo {year} {2014}),\ 10.1017/jfm.2014.685}\BibitemShut {NoStop}%
\bibitem [{\citenamefont {Bühler}, \citenamefont {Kuang},\ and\ \citenamefont {Tabak}(2017)}]{Buhler2017361}%
  \BibitemOpen
  \bibfield  {author} {\bibinfo {author} {\bibfnamefont {O.}~\bibnamefont {Bühler}}, \bibinfo {author} {\bibfnamefont {M.}~\bibnamefont {Kuang}},\ and\ \bibinfo {author} {\bibfnamefont {E.}~\bibnamefont {Tabak}},\ }\bibfield  {title} {\enquote {\bibinfo {title} {Anisotropic {H}elmholtz and wave-vortex decomposition of one-dimensional spectra},}\ }\href {https://doi.org/10.1017/jfm.2017.57} {\bibfield  {journal} {\bibinfo  {journal} {J. Fluid Mech.}\ }\textbf {\bibinfo {volume} {815}},\ \bibinfo {pages} {361--387} (\bibinfo {year} {2017})}\BibitemShut {NoStop}%
\bibitem [{\citenamefont {Schoder}, \citenamefont {Roppert},\ and\ \citenamefont {Kaltenbacher}(2020)}]{Schoder20203019}%
  \BibitemOpen
  \bibfield  {author} {\bibinfo {author} {\bibfnamefont {S.}~\bibnamefont {Schoder}}, \bibinfo {author} {\bibfnamefont {K.}~\bibnamefont {Roppert}},\ and\ \bibinfo {author} {\bibfnamefont {M.}~\bibnamefont {Kaltenbacher}},\ }\bibfield  {title} {\enquote {\bibinfo {title} {Postprocessing of direct aeroacoustic simulations using {H}elmholtz decomposition},}\ }\href {https://doi.org/10.2514/1.J058836} {\bibfield  {journal} {\bibinfo  {journal} {AIAA J.}\ }\textbf {\bibinfo {volume} {58}},\ \bibinfo {pages} {3019--3027} (\bibinfo {year} {2020})}\BibitemShut {NoStop}%
\bibitem [{\citenamefont {Caltagirone}(2021)}]{Caltagirone2021}%
  \BibitemOpen
  \bibfield  {author} {\bibinfo {author} {\bibfnamefont {J.-P.}\ \bibnamefont {Caltagirone}},\ }\bibfield  {title} {\enquote {\bibinfo {title} {On a reformulation of {N}avier--{S}tokes equations based on {H}elmholtz--{H}odge decomposition},}\ }\href {https://doi.org/10.1063/5.0053412} {\bibfield  {journal} {\bibinfo  {journal} {Phys. Fluids}\ }\textbf {\bibinfo {volume} {33}} (\bibinfo {year} {2021}),\ 10.1063/5.0053412}\BibitemShut {NoStop}%
\bibitem [{\citenamefont {Aharonov}\ and\ \citenamefont {Bohm}(1959)}]{Aharonov1959485}%
  \BibitemOpen
  \bibfield  {author} {\bibinfo {author} {\bibfnamefont {Y.}~\bibnamefont {Aharonov}}\ and\ \bibinfo {author} {\bibfnamefont {D.}~\bibnamefont {Bohm}},\ }\bibfield  {title} {\enquote {\bibinfo {title} {Significance of electromagnetic potentials in the quantum theory},}\ }\href {https://doi.org/10.1103/PhysRev.115.485} {\bibfield  {journal} {\bibinfo  {journal} {Phys. Rev.}\ }\textbf {\bibinfo {volume} {115}},\ \bibinfo {pages} {485--491} (\bibinfo {year} {1959})}\BibitemShut {NoStop}%
\bibitem [{\citenamefont {Konopinski}(1978)}]{Konopinski1978499}%
  \BibitemOpen
  \bibfield  {author} {\bibinfo {author} {\bibfnamefont {E.}~\bibnamefont {Konopinski}},\ }\bibfield  {title} {\enquote {\bibinfo {title} {What the electromagnetic vector potential describes},}\ }\href {https://doi.org/10.1119/1.11298} {\bibfield  {journal} {\bibinfo  {journal} {Am. J. Phys.}\ }\textbf {\bibinfo {volume} {46}},\ \bibinfo {pages} {499--502} (\bibinfo {year} {1978})}\BibitemShut {NoStop}%
\bibitem [{\citenamefont {Haber}\ \emph {et~al.}(2000)\citenamefont {Haber}, \citenamefont {Ascher}, \citenamefont {Aruliah},\ and\ \citenamefont {Oldenburg}}]{Haber2000150}%
  \BibitemOpen
  \bibfield  {author} {\bibinfo {author} {\bibfnamefont {E.}~\bibnamefont {Haber}}, \bibinfo {author} {\bibfnamefont {U.}~\bibnamefont {Ascher}}, \bibinfo {author} {\bibfnamefont {D.}~\bibnamefont {Aruliah}},\ and\ \bibinfo {author} {\bibfnamefont {D.}~\bibnamefont {Oldenburg}},\ }\bibfield  {title} {\enquote {\bibinfo {title} {Fast simulation of 3d electromagnetic problems using potentials},}\ }\href {https://doi.org/10.1006/jcph.2000.6545} {\bibfield  {journal} {\bibinfo  {journal} {J. Comput. Phys.}\ }\textbf {\bibinfo {volume} {163}},\ \bibinfo {pages} {150--171} (\bibinfo {year} {2000})}\BibitemShut {NoStop}%
\bibitem [{\citenamefont {Hehl}\ and\ \citenamefont {Obukhov}(2003)}]{hehl2003foundations}%
  \BibitemOpen
  \bibfield  {author} {\bibinfo {author} {\bibfnamefont {F.~W.}\ \bibnamefont {Hehl}}\ and\ \bibinfo {author} {\bibfnamefont {Y.~N.}\ \bibnamefont {Obukhov}},\ }\href@noop {} {\emph {\bibinfo {title} {Foundations of classical electrodynamics: Charge, flux, and metric}}},\ Vol.~\bibinfo {volume} {33}\ (\bibinfo  {publisher} {Springer {S}cience \& {B}usiness {M}edia},\ \bibinfo {year} {2003})\BibitemShut {NoStop}%
\bibitem [{\citenamefont {Weiss}(2013)}]{Weiss201340}%
  \BibitemOpen
  \bibfield  {author} {\bibinfo {author} {\bibfnamefont {C.}~\bibnamefont {Weiss}},\ }\bibfield  {title} {\enquote {\bibinfo {title} {Project {APhiD}: A {L}orenz-gauged ${A}-\phi$ decomposition for parallelized computation of ultra-broadband electromagnetic induction in a fully heterogeneous earth},}\ }\href {https://doi.org/10.1016/j.cageo.2013.05.002} {\bibfield  {journal} {\bibinfo  {journal} {Comput. Geosci.}\ }\textbf {\bibinfo {volume} {58}},\ \bibinfo {pages} {40--52} (\bibinfo {year} {2013})}\BibitemShut {NoStop}%
\bibitem [{\citenamefont {Du}\ \emph {et~al.}(2017)\citenamefont {Du}, \citenamefont {Guo}, \citenamefont {Zhao}, \citenamefont {Gong}, \citenamefont {Wang},\ and\ \citenamefont {Li}}]{Du2017S111}%
  \BibitemOpen
  \bibfield  {author} {\bibinfo {author} {\bibfnamefont {Q.}~\bibnamefont {Du}}, \bibinfo {author} {\bibfnamefont {C.}~\bibnamefont {Guo}}, \bibinfo {author} {\bibfnamefont {Q.}~\bibnamefont {Zhao}}, \bibinfo {author} {\bibfnamefont {X.}~\bibnamefont {Gong}}, \bibinfo {author} {\bibfnamefont {C.}~\bibnamefont {Wang}},\ and\ \bibinfo {author} {\bibfnamefont {X.-Y.}\ \bibnamefont {Li}},\ }\bibfield  {title} {\enquote {\bibinfo {title} {Vector-based elastic reverse time migration based on scalar imaging condition},}\ }\href {https://doi.org/10.1190/GEO2016-0146.1} {\bibfield  {journal} {\bibinfo  {journal} {Geophysics}\ }\textbf {\bibinfo {volume} {82}},\ \bibinfo {pages} {S111--S127} (\bibinfo {year} {2017})}\BibitemShut {NoStop}%
\bibitem [{\citenamefont {Shi}, \citenamefont {Zhang},\ and\ \citenamefont {Wang}(2019)}]{Shi2019509}%
  \BibitemOpen
  \bibfield  {author} {\bibinfo {author} {\bibfnamefont {Y.}~\bibnamefont {Shi}}, \bibinfo {author} {\bibfnamefont {W.}~\bibnamefont {Zhang}},\ and\ \bibinfo {author} {\bibfnamefont {Y.}~\bibnamefont {Wang}},\ }\bibfield  {title} {\enquote {\bibinfo {title} {Seismic elastic {RTM} with vector-wavefield decomposition},}\ }\href {https://doi.org/10.1093/jge/gxz023} {\bibfield  {journal} {\bibinfo  {journal} {J. Geophys. Eng.}\ }\textbf {\bibinfo {volume} {16}},\ \bibinfo {pages} {509--524} (\bibinfo {year} {2019})}\BibitemShut {NoStop}%
\bibitem [{\citenamefont {Paganin}\ and\ \citenamefont {Nugent}(1998)}]{Paganin19982586}%
  \BibitemOpen
  \bibfield  {author} {\bibinfo {author} {\bibfnamefont {D.}~\bibnamefont {Paganin}}\ and\ \bibinfo {author} {\bibfnamefont {K.}~\bibnamefont {Nugent}},\ }\bibfield  {title} {\enquote {\bibinfo {title} {Noninterferometric phase imaging with partially coherent light},}\ }\href {https://doi.org/10.1103/PhysRevLett.80.2586} {\bibfield  {journal} {\bibinfo  {journal} {Phy. Rev. Lett.}\ }\textbf {\bibinfo {volume} {80}},\ \bibinfo {pages} {2586--2589} (\bibinfo {year} {1998})}\BibitemShut {NoStop}%
\bibitem [{\citenamefont {Park}\ and\ \citenamefont {Maniatty}(2006)}]{Park20063697}%
  \BibitemOpen
  \bibfield  {author} {\bibinfo {author} {\bibfnamefont {E.}~\bibnamefont {Park}}\ and\ \bibinfo {author} {\bibfnamefont {A.}~\bibnamefont {Maniatty}},\ }\bibfield  {title} {\enquote {\bibinfo {title} {Shear modulus reconstruction in dynamic elastography: Time harmonic case},}\ }\href {https://doi.org/10.1088/0031-9155/51/15/007} {\bibfield  {journal} {\bibinfo  {journal} {Phys. Med. Biol.}\ }\textbf {\bibinfo {volume} {51}},\ \bibinfo {pages} {3697--3721} (\bibinfo {year} {2006})}\BibitemShut {NoStop}%
\bibitem [{\citenamefont {Kohlberger}, \citenamefont {M{\'e}min},\ and\ \citenamefont {Schnörr}(2003)}]{Kohlberger2003432}%
  \BibitemOpen
  \bibfield  {author} {\bibinfo {author} {\bibfnamefont {T.}~\bibnamefont {Kohlberger}}, \bibinfo {author} {\bibfnamefont {E.}~\bibnamefont {M{\'e}min}},\ and\ \bibinfo {author} {\bibfnamefont {C.}~\bibnamefont {Schnörr}},\ }\bibfield  {title} {\enquote {\bibinfo {title} {Variational dense motion estimation using the {H}elmholtz decomposition},}\ }\href {https://doi.org/10.1007/3-540-44935-3_30} {\bibfield  {journal} {\bibinfo  {journal} {Lect. Notes Comput. Sci.}\ }\textbf {\bibinfo {volume} {2695}},\ \bibinfo {pages} {432--448} (\bibinfo {year} {2003})}\BibitemShut {NoStop}%
\bibitem [{\citenamefont {Guo}, \citenamefont {Mandal},\ and\ \citenamefont {Li}(2005)}]{Guo2005493}%
  \BibitemOpen
  \bibfield  {author} {\bibinfo {author} {\bibfnamefont {Q.}~\bibnamefont {Guo}}, \bibinfo {author} {\bibfnamefont {M.}~\bibnamefont {Mandal}},\ and\ \bibinfo {author} {\bibfnamefont {M.}~\bibnamefont {Li}},\ }\bibfield  {title} {\enquote {\bibinfo {title} {Efficient {H}odge-{H}elmholtz decomposition of motion fields},}\ }\href {https://doi.org/10.1016/j.patrec.2004.08.008} {\bibfield  {journal} {\bibinfo  {journal} {Pattern Recognit. Lett.}\ }\textbf {\bibinfo {volume} {26}},\ \bibinfo {pages} {493--501} (\bibinfo {year} {2005})}\BibitemShut {NoStop}%
\bibitem [{\citenamefont {Cuzol}, \citenamefont {Hellier},\ and\ \citenamefont {M{\'e}min}(2007)}]{Cuzol2007329}%
  \BibitemOpen
  \bibfield  {author} {\bibinfo {author} {\bibfnamefont {A.}~\bibnamefont {Cuzol}}, \bibinfo {author} {\bibfnamefont {P.}~\bibnamefont {Hellier}},\ and\ \bibinfo {author} {\bibfnamefont {E.}~\bibnamefont {M{\'e}min}},\ }\bibfield  {title} {\enquote {\bibinfo {title} {A low dimensional fluid motion estimator},}\ }\href {https://doi.org/10.1007/s11263-007-0037-0} {\bibfield  {journal} {\bibinfo  {journal} {Int. J. Comput. Vis.}\ }\textbf {\bibinfo {volume} {75}},\ \bibinfo {pages} {329--349} (\bibinfo {year} {2007})}\BibitemShut {NoStop}%
\bibitem [{\citenamefont {Simader}\ and\ \citenamefont {Sohr}(1992)}]{simader92}%
  \BibitemOpen
  \bibfield  {author} {\bibinfo {author} {\bibfnamefont {C.~G.}\ \bibnamefont {Simader}}\ and\ \bibinfo {author} {\bibfnamefont {H.}~\bibnamefont {Sohr}},\ }\enquote {\bibinfo {title} {A new approach to the {H}elmholtz decomposition and the {N}eumann problem in {L}\textsuperscript{q}-spaces for bounded and exterior domains},}\ in\ \href {https://doi.org/10.1142/9789814503594_0001} {\emph {\bibinfo {booktitle} {Mathematical Problems Relating to the Navier-Stokes Equations}}},\ \bibinfo {series} {Series on Advances in Mathematics for Applied Sciences}, Vol.~\bibinfo {volume} {11}\ (\bibinfo  {publisher} {World {S}cientific {P}ublishing},\ \bibinfo {year} {1992})\ pp.\ \bibinfo {pages} {1--35}\BibitemShut {NoStop}%
\bibitem [{\citenamefont {Farwig}, \citenamefont {Kozono},\ and\ \citenamefont {Sohr}(2007)}]{Farwig2007239}%
  \BibitemOpen
  \bibfield  {author} {\bibinfo {author} {\bibfnamefont {R.}~\bibnamefont {Farwig}}, \bibinfo {author} {\bibfnamefont {H.}~\bibnamefont {Kozono}},\ and\ \bibinfo {author} {\bibfnamefont {H.}~\bibnamefont {Sohr}},\ }\bibfield  {title} {\enquote {\bibinfo {title} {On the {H}elmholtz decomposition in general unbounded domains},}\ }\href {https://doi.org/10.1007/s00013-006-1910-8} {\bibfield  {journal} {\bibinfo  {journal} {Arch. Math.}\ }\textbf {\bibinfo {volume} {88}},\ \bibinfo {pages} {239--248} (\bibinfo {year} {2007})}\BibitemShut {NoStop}%
\bibitem [{\citenamefont {Sohr}(2012)}]{sohr2012navier}%
  \BibitemOpen
  \bibfield  {author} {\bibinfo {author} {\bibfnamefont {H.}~\bibnamefont {Sohr}},\ }\href@noop {} {\emph {\bibinfo {title} {The {N}avier--{S}tokes equations: An elementary functional analytic approach}}}\ (\bibinfo  {publisher} {Springer {S}cience \& {B}usiness {M}edia},\ \bibinfo {year} {2012})\BibitemShut {NoStop}%
\bibitem [{\citenamefont {Moses}(1971)}]{moses1971}%
  \BibitemOpen
  \bibfield  {author} {\bibinfo {author} {\bibfnamefont {H.~E.}\ \bibnamefont {Moses}},\ }\bibfield  {title} {\enquote {\bibinfo {title} {Eigenfunctions of the curl operator, rotationally invariant helmholtz theorem, and applications to electromagnetic theory and fluid mechanics},}\ }\href {https://doi.org/10.1137/0121015} {\bibfield  {journal} {\bibinfo  {journal} {SIAM J. Appl. Math.}\ }\textbf {\bibinfo {volume} {21}},\ \bibinfo {pages} {114--144} (\bibinfo {year} {1971})}\BibitemShut {NoStop}%
\bibitem [{\citenamefont {Schweizer}(2018)}]{Schweizer2018}%
  \BibitemOpen
  \bibfield  {author} {\bibinfo {author} {\bibfnamefont {B.}~\bibnamefont {Schweizer}},\ }\enquote {\bibinfo {title} {On {Friedrichs} inequality, {H}elmholtz decomposition, vector potentials, and the div-curl lemma},}\ in\ \href {https://doi.org/10.1007/978-3-319-75940-1_4} {\emph {\bibinfo {booktitle} {Trends in Applications of Mathematics to Mechanics}}},\ \bibinfo {editor} {edited by\ \bibinfo {editor} {\bibfnamefont {E.}~\bibnamefont {Rocca}}, \bibinfo {editor} {\bibfnamefont {U.}~\bibnamefont {Stefanelli}}, \bibinfo {editor} {\bibfnamefont {L.}~\bibnamefont {Truskinovsky}},\ and\ \bibinfo {editor} {\bibfnamefont {A.}~\bibnamefont {Visintin}}}\ (\bibinfo  {publisher} {Springer {I}nternational {P}ublishing},\ \bibinfo {address} {Cham},\ \bibinfo {year} {2018})\ pp.\ \bibinfo {pages} {65--79}\BibitemShut {NoStop}%
\bibitem [{\citenamefont {Giga}\ and\ \citenamefont {Gu}(2022)}]{Giga2022}%
  \BibitemOpen
  \bibfield  {author} {\bibinfo {author} {\bibfnamefont {Y.}~\bibnamefont {Giga}}\ and\ \bibinfo {author} {\bibfnamefont {Z.}~\bibnamefont {Gu}},\ }\bibfield  {title} {\enquote {\bibinfo {title} {The {H}elmholtz decomposition of a space of vector fields with bounded mean oscillation in a bounded domain},}\ }\href {https://doi.org/10.1007/s00208-022-02410-y} {\bibfield  {journal} {\bibinfo  {journal} {Math. Ann.}\ } (\bibinfo {year} {2022}),\ 10.1007/s00208-022-02410-y}\BibitemShut {NoStop}%
\bibitem [{\citenamefont {Glötzl}\ and\ \citenamefont {Richters}(2023)}]{Glotzl2023}%
  \BibitemOpen
  \bibfield  {author} {\bibinfo {author} {\bibfnamefont {E.}~\bibnamefont {Glötzl}}\ and\ \bibinfo {author} {\bibfnamefont {O.}~\bibnamefont {Richters}},\ }\bibfield  {title} {\enquote {\bibinfo {title} {{H}elmholtz decomposition and potential functions for n-dimensional analytic vector fields},}\ }\href {https://doi.org/10.1016/j.jmaa.2023.127138} {\bibfield  {journal} {\bibinfo  {journal} {J. Math. Anal. Appl.}\ }\textbf {\bibinfo {volume} {525}} (\bibinfo {year} {2023}),\ 10.1016/j.jmaa.2023.127138}\BibitemShut {NoStop}%
\bibitem [{\citenamefont {Bhatia}\ \emph {et~al.}(2013)\citenamefont {Bhatia}, \citenamefont {Norgard}, \citenamefont {Pascucci},\ and\ \citenamefont {Bremer}}]{Bhatia20131386}%
  \BibitemOpen
  \bibfield  {author} {\bibinfo {author} {\bibfnamefont {H.}~\bibnamefont {Bhatia}}, \bibinfo {author} {\bibfnamefont {G.}~\bibnamefont {Norgard}}, \bibinfo {author} {\bibfnamefont {V.}~\bibnamefont {Pascucci}},\ and\ \bibinfo {author} {\bibfnamefont {P.-T.}\ \bibnamefont {Bremer}},\ }\bibfield  {title} {\enquote {\bibinfo {title} {The {H}elmholtz-hodge decomposition - a survey},}\ }\href {https://doi.org/10.1109/TVCG.2012.316} {\bibfield  {journal} {\bibinfo  {journal} {IEEE Trans. Vis. Comput. Graph.}\ }\textbf {\bibinfo {volume} {19}},\ \bibinfo {pages} {1386--1404} (\bibinfo {year} {2013})}\BibitemShut {NoStop}%
\bibitem [{\citenamefont {Horn}\ and\ \citenamefont {Johnson}(2012)}]{horn2012matrix}%
  \BibitemOpen
  \bibfield  {author} {\bibinfo {author} {\bibfnamefont {R.~A.}\ \bibnamefont {Horn}}\ and\ \bibinfo {author} {\bibfnamefont {C.~R.}\ \bibnamefont {Johnson}},\ }\href@noop {} {\emph {\bibinfo {title} {Matrix analysis}}}\ (\bibinfo  {publisher} {Cambridge {U}niversity press},\ \bibinfo {year} {2012})\BibitemShut {NoStop}%
\bibitem [{Note2()}]{Note2}%
  \BibitemOpen
  \bibinfo {note} {Note that, by the uniqueness of the Cholesky decomposition for positive definite matrices,\cite [see p. 441]{horn2012matrix} the polar decomposition $J=Qh$ and the QR factorization\cite [see p. 449]{horn2012matrix} of $J$ coincide, which implies that $h^{-1}$ is generally a lower triangular matrix.}\BibitemShut {Stop}%
\bibitem [{\citenamefont {Pedergnana}\ and\ \citenamefont {Noiray}(2022)}]{Pedergnana2022}%
  \BibitemOpen
  \bibfield  {author} {\bibinfo {author} {\bibfnamefont {T.}~\bibnamefont {Pedergnana}}\ and\ \bibinfo {author} {\bibfnamefont {N.}~\bibnamefont {Noiray}},\ }\bibfield  {title} {\enquote {\bibinfo {title} {Exact potentials in multivariate {L}angevin equations},}\ }\href {https://doi.org/10.1063/5.0124031} {\bibfield  {journal} {\bibinfo  {journal} {Chaos}\ }\textbf {\bibinfo {volume} {32}} (\bibinfo {year} {2022}),\ 10.1063/5.0124031}\BibitemShut {NoStop}%
\bibitem [{\citenamefont {Guillemin}\ and\ \citenamefont {Pollack}(2010)}]{guillemin2010differential}%
  \BibitemOpen
  \bibfield  {author} {\bibinfo {author} {\bibfnamefont {V.}~\bibnamefont {Guillemin}}\ and\ \bibinfo {author} {\bibfnamefont {A.}~\bibnamefont {Pollack}},\ }\href@noop {} {\emph {\bibinfo {title} {Differential topology}}},\ Vol.\ \bibinfo {volume} {370}\ (\bibinfo  {publisher} {American Mathematical Society},\ \bibinfo {year} {2010})\BibitemShut {NoStop}%
\bibitem [{\citenamefont {Gonz\'{a}lez-Gasc\'{o}n}(1986)}]{Gonzalez-Gascon198661}%
  \BibitemOpen
  \bibfield  {author} {\bibinfo {author} {\bibfnamefont {F.}~\bibnamefont {Gonz\'{a}lez-Gasc\'{o}n}},\ }\bibfield  {title} {\enquote {\bibinfo {title} {Note on paper of {A}ndrey concerning non-{H}amiltonian systems},}\ }\href {https://doi.org/10.1016/0375-9601(86)90478-0} {\bibfield  {journal} {\bibinfo  {journal} {Phys. Lett. A}\ }\textbf {\bibinfo {volume} {114}},\ \bibinfo {pages} {61 – 62} (\bibinfo {year} {1986})}\BibitemShut {NoStop}%
\bibitem [{\citenamefont {Nutku}(1990{\natexlab{a}})}]{Nutku199027}%
  \BibitemOpen
  \bibfield  {author} {\bibinfo {author} {\bibfnamefont {Y.}~\bibnamefont {Nutku}},\ }\bibfield  {title} {\enquote {\bibinfo {title} {Hamiltonian structure of the {L}otka-{V}olterra equations},}\ }\href {https://doi.org/10.1016/0375-9601(90)90270-X} {\bibfield  {journal} {\bibinfo  {journal} {Phys. Lett. A}\ }\textbf {\bibinfo {volume} {145}},\ \bibinfo {pages} {27 – 28} (\bibinfo {year} {1990}{\natexlab{a}})}\BibitemShut {NoStop}%
\bibitem [{Note3()}]{Note3}%
  \BibitemOpen
  \bibinfo {note} {An ``orbit'' is {a curve} in the (frozen) phase space of a steady dynamical system along which a trajectory runs.}\BibitemShut {Stop}%
\bibitem [{\citenamefont {del Pino}, \citenamefont {Košata},\ and\ \citenamefont {Zilberberg}(2023)}]{delpino2023limit}%
  \BibitemOpen
  \bibfield  {author} {\bibinfo {author} {\bibfnamefont {J.}~\bibnamefont {del Pino}}, \bibinfo {author} {\bibfnamefont {J.}~\bibnamefont {Košata}},\ and\ \bibinfo {author} {\bibfnamefont {O.}~\bibnamefont {Zilberberg}},\ }\href@noop {} {\enquote {\bibinfo {title} {Limit cycles as stationary states of an extended harmonic balance ansatz},}\ } (\bibinfo {year} {2023}),\ \Eprint {https://arxiv.org/abs/2308.06092} {arXiv:2308.06092 [nlin.AO]} \BibitemShut {NoStop}%
\bibitem [{\citenamefont {Parker}\ and\ \citenamefont {Chua}(2012)}]{parker2012practical}%
  \BibitemOpen
  \bibfield  {author} {\bibinfo {author} {\bibfnamefont {T.~S.}\ \bibnamefont {Parker}}\ and\ \bibinfo {author} {\bibfnamefont {L.}~\bibnamefont {Chua}},\ }\href@noop {} {\emph {\bibinfo {title} {Practical numerical algorithms for chaotic systems}}}\ (\bibinfo  {publisher} {Springer Science \& Business Media},\ \bibinfo {year} {2012})\BibitemShut {NoStop}%
\bibitem [{\citenamefont {Han}(2006)}]{han2006bifurcation}%
  \BibitemOpen
  \bibfield  {author} {\bibinfo {author} {\bibfnamefont {M.}~\bibnamefont {Han}},\ }\bibfield  {title} {\enquote {\bibinfo {title} {Bifurcation theory of limit cycles of planar systems},}\ }in\ \href@noop {} {\emph {\bibinfo {booktitle} {Handbook of Differential Equations: Ordinary Differential Equations}}},\ Vol.~\bibinfo {volume} {3}\ (\bibinfo  {publisher} {Elsevier},\ \bibinfo {year} {2006})\ pp.\ \bibinfo {pages} {341--433}\BibitemShut {NoStop}%
\bibitem [{\citenamefont {Hilbert}(1900)}]{hilbert1900mathematische}%
  \BibitemOpen
  \bibfield  {author} {\bibinfo {author} {\bibfnamefont {D.}~\bibnamefont {Hilbert}},\ }\bibfield  {title} {\enquote {\bibinfo {title} {Mathematische {P}robleme},}\ }\href@noop {} {\bibfield  {journal} {\bibinfo  {journal} {Nachrichten von der {Kö}niglichen {G}esellschaft der {W}issenschaften zu {Gö}ttingen}\ } (\bibinfo {year} {1900})}\BibitemShut {NoStop}%
\bibitem [{\citenamefont {Ilyashenko}(2002)}]{Ilyashenko2002301}%
  \BibitemOpen
  \bibfield  {author} {\bibinfo {author} {\bibfnamefont {Y.}~\bibnamefont {Ilyashenko}},\ }\bibfield  {title} {\enquote {\bibinfo {title} {Centennial history of hilbert's 16th problem},}\ }\href {https://doi.org/10.1090/S0273-0979-02-00946-1} {\bibfield  {journal} {\bibinfo  {journal} {Bull. Am. Math. Soc.}\ }\textbf {\bibinfo {volume} {39}},\ \bibinfo {pages} {301 – 354} (\bibinfo {year} {2002})}\BibitemShut {NoStop}%
\bibitem [{\citenamefont {Busenberg}\ and\ \citenamefont {Van~Den}(1993)}]{Busenberg1993463}%
  \BibitemOpen
  \bibfield  {author} {\bibinfo {author} {\bibfnamefont {S.}~\bibnamefont {Busenberg}}\ and\ \bibinfo {author} {\bibfnamefont {P.~D.}\ \bibnamefont {Van~Den}},\ }\bibfield  {title} {\enquote {\bibinfo {title} {A method for proving the non-existence of limit cycles},}\ }\href {https://doi.org/10.1006/jmaa.1993.1037} {\bibfield  {journal} {\bibinfo  {journal} {J. Math. Anal. Appl.}\ }\textbf {\bibinfo {volume} {172}},\ \bibinfo {pages} {463 – 479} (\bibinfo {year} {1993})}\BibitemShut {NoStop}%
\bibitem [{Note4()}]{Note4}%
  \BibitemOpen
  \bibinfo {note} {The {limiting} case of a closed orbit which is a fixed point is excluded from the discussion here.}\BibitemShut {Stop}%
\bibitem [{\citenamefont {Jost}(2011)}]{jost2008riemannian}%
  \BibitemOpen
  \bibfield  {author} {\bibinfo {author} {\bibfnamefont {J.}~\bibnamefont {Jost}},\ }\href@noop {} {\emph {\bibinfo {title} {Riemannian geometry and geometric analysis}}},\ Vol.\ \bibinfo {volume} {42005}\ (\bibinfo  {publisher} {Springer},\ \bibinfo {year} {2011})\BibitemShut {NoStop}%
\bibitem [{\citenamefont {Shi}(1980)}]{shi1980concrete}%
  \BibitemOpen
  \bibfield  {author} {\bibinfo {author} {\bibfnamefont {S.}~\bibnamefont {Shi}},\ }\bibfield  {title} {\enquote {\bibinfo {title} {Concrete example of the existence of 4 limit-cycles for plane quadratic systems},}\ }\href@noop {} {\bibfield  {journal} {\bibinfo  {journal} {Sci. Sin.}\ }\textbf {\bibinfo {volume} {23}},\ \bibinfo {pages} {153--158} (\bibinfo {year} {1980})}\BibitemShut {NoStop}%
\bibitem [{\citenamefont {Galias}\ and\ \citenamefont {Tucker}(2022)}]{Galias2022}%
  \BibitemOpen
  \bibfield  {author} {\bibinfo {author} {\bibfnamefont {Z.}~\bibnamefont {Galias}}\ and\ \bibinfo {author} {\bibfnamefont {W.}~\bibnamefont {Tucker}},\ }\bibfield  {title} {\enquote {\bibinfo {title} {The {S}ongling system has exactly four limit cycles},}\ }\href {https://doi.org/10.1016/j.amc.2021.126691} {\bibfield  {journal} {\bibinfo  {journal} {Appl. Math. Comput.}\ }\textbf {\bibinfo {volume} {415}} (\bibinfo {year} {2022}),\ 10.1016/j.amc.2021.126691}\BibitemShut {NoStop}%
\bibitem [{\citenamefont {S{\"a}rkk{\"a}}\ and\ \citenamefont {Solin}(2019)}]{sarkka2019applied}%
  \BibitemOpen
  \bibfield  {author} {\bibinfo {author} {\bibfnamefont {S.}~\bibnamefont {S{\"a}rkk{\"a}}}\ and\ \bibinfo {author} {\bibfnamefont {A.}~\bibnamefont {Solin}},\ }\href@noop {} {\emph {\bibinfo {title} {Applied stochastic differential equations}}},\ Vol.~\bibinfo {volume} {10}\ (\bibinfo  {publisher} {Cambridge University Press},\ \bibinfo {year} {2019})\BibitemShut {NoStop}%
\bibitem [{\citenamefont {Abraham}, \citenamefont {Marsden},\ and\ \citenamefont {Ratiu}(2012)}]{abraham2012manifolds}%
  \BibitemOpen
  \bibfield  {author} {\bibinfo {author} {\bibfnamefont {R.}~\bibnamefont {Abraham}}, \bibinfo {author} {\bibfnamefont {J.~E.}\ \bibnamefont {Marsden}},\ and\ \bibinfo {author} {\bibfnamefont {T.}~\bibnamefont {Ratiu}},\ }\href@noop {} {\emph {\bibinfo {title} {Manifolds, tensor analysis, and applications}}},\ Vol.~\bibinfo {volume} {75}\ (\bibinfo  {publisher} {Springer Science \& Business Media},\ \bibinfo {year} {2012})\BibitemShut {NoStop}%
\bibitem [{\citenamefont {Olver}(1993)}]{olver1993applications}%
  \BibitemOpen
  \bibfield  {author} {\bibinfo {author} {\bibfnamefont {P.~J.}\ \bibnamefont {Olver}},\ }\href@noop {} {\emph {\bibinfo {title} {Applications of {L}ie groups to differential equations}}},\ Vol.\ \bibinfo {volume} {107}\ (\bibinfo  {publisher} {Springer Science \& Business Media},\ \bibinfo {year} {1993})\BibitemShut {NoStop}%
\bibitem [{\citenamefont {Quispel}\ and\ \citenamefont {Capel}(1996)}]{Quispel1996223}%
  \BibitemOpen
  \bibfield  {author} {\bibinfo {author} {\bibfnamefont {G.}~\bibnamefont {Quispel}}\ and\ \bibinfo {author} {\bibfnamefont {H.}~\bibnamefont {Capel}},\ }\bibfield  {title} {\enquote {\bibinfo {title} {Solving odes numerically while preserving a first integral},}\ }\href {https://doi.org/10.1016/0375-9601(96)00403-3} {\bibfield  {journal} {\bibinfo  {journal} {Phys. Lett. A}\ }\textbf {\bibinfo {volume} {218}},\ \bibinfo {pages} {223 – 228} (\bibinfo {year} {1996})}\BibitemShut {NoStop}%
\bibitem [{\citenamefont {McLachlan}, \citenamefont {Quispel},\ and\ \citenamefont {Robidoux}(1998)}]{McLachlan19982399}%
  \BibitemOpen
  \bibfield  {author} {\bibinfo {author} {\bibfnamefont {R.~I.}\ \bibnamefont {McLachlan}}, \bibinfo {author} {\bibfnamefont {G.}~\bibnamefont {Quispel}},\ and\ \bibinfo {author} {\bibfnamefont {N.}~\bibnamefont {Robidoux}},\ }\bibfield  {title} {\enquote {\bibinfo {title} {Unified approach to {H}amiltonian systems, {P}oisson systems, gradient systems, and systems with {L}yapunov functions or first integrals},}\ }\href {https://doi.org/10.1103/PhysRevLett.81.2399} {\bibfield  {journal} {\bibinfo  {journal} {Phys. Rev. Lett.}\ }\textbf {\bibinfo {volume} {81}},\ \bibinfo {pages} {2399 – 2403} (\bibinfo {year} {1998})}\BibitemShut {NoStop}%
\bibitem [{\citenamefont {McLachlan}, \citenamefont {Quispel},\ and\ \citenamefont {Robidoux}(1999)}]{McLachlan19991021}%
  \BibitemOpen
  \bibfield  {author} {\bibinfo {author} {\bibfnamefont {R.~I.}\ \bibnamefont {McLachlan}}, \bibinfo {author} {\bibfnamefont {G.}~\bibnamefont {Quispel}},\ and\ \bibinfo {author} {\bibfnamefont {N.}~\bibnamefont {Robidoux}},\ }\bibfield  {title} {\enquote {\bibinfo {title} {Geometric integration using discrete gradients},}\ }\href {https://doi.org/10.1098/rsta.1999.0363} {\bibfield  {journal} {\bibinfo  {journal} {Philos. Trans. R. Soc. A}\ }\textbf {\bibinfo {volume} {357}},\ \bibinfo {pages} {1021 – 1045} (\bibinfo {year} {1999})}\BibitemShut {NoStop}%
\bibitem [{\citenamefont {Bárta}, \citenamefont {Chill},\ and\ \citenamefont {Fašangová}(2012)}]{Barta201257}%
  \BibitemOpen
  \bibfield  {author} {\bibinfo {author} {\bibfnamefont {T.}~\bibnamefont {Bárta}}, \bibinfo {author} {\bibfnamefont {R.}~\bibnamefont {Chill}},\ and\ \bibinfo {author} {\bibfnamefont {E.}~\bibnamefont {Fašangová}},\ }\bibfield  {title} {\enquote {\bibinfo {title} {Every ordinary differential equation with a strict {L}yapunov function is a gradient system},}\ }\href {https://doi.org/10.1007/s00605-011-0322-4} {\bibfield  {journal} {\bibinfo  {journal} {Monatsh. Math.}\ }\textbf {\bibinfo {volume} {166}},\ \bibinfo {pages} {57 – 72} (\bibinfo {year} {2012})}\BibitemShut {NoStop}%
\bibitem [{Note5()}]{Note5}%
  \BibitemOpen
  \bibinfo {note} {Unless $u$ vanishes identically, in which case there is nothing to prove.}\BibitemShut {Stop}%
\bibitem [{\citenamefont {Kermack}\ and\ \citenamefont {McKendrick}(1991)}]{Kermack199133}%
  \BibitemOpen
  \bibfield  {author} {\bibinfo {author} {\bibfnamefont {W.}~\bibnamefont {Kermack}}\ and\ \bibinfo {author} {\bibfnamefont {A.}~\bibnamefont {McKendrick}},\ }\bibfield  {title} {\enquote {\bibinfo {title} {Contributions to the mathematical theory of epidemics-{I}},}\ }\href {https://doi.org/10.1007/BF02464423} {\bibfield  {journal} {\bibinfo  {journal} {Bull. Math. Biol.}\ }\textbf {\bibinfo {volume} {53}},\ \bibinfo {pages} {33 – 55} (\bibinfo {year} {1991})}\BibitemShut {NoStop}%
\bibitem [{\citenamefont {Nutku}(1990{\natexlab{b}})}]{Nutku1990L1145}%
  \BibitemOpen
  \bibfield  {author} {\bibinfo {author} {\bibfnamefont {Y.}~\bibnamefont {Nutku}},\ }\bibfield  {title} {\enquote {\bibinfo {title} {Bi-{H}amiltonian structure of the {K}ermack-{M}c{K}endrick model for epidemics},}\ }\href {https://doi.org/10.1088/0305-4470/23/21/013} {\bibfield  {journal} {\bibinfo  {journal} {Journal of Physics A: Mathematical and General}\ }\textbf {\bibinfo {volume} {23}},\ \bibinfo {pages} {L1145–L1146} (\bibinfo {year} {1990}{\natexlab{b}})}\BibitemShut {NoStop}%
\bibitem [{\citenamefont {Ballesteros}, \citenamefont {Blasco},\ and\ \citenamefont {Gutierrez-Sagredo}(2020)}]{Ballesteros2020}%
  \BibitemOpen
  \bibfield  {author} {\bibinfo {author} {\bibfnamefont {A.}~\bibnamefont {Ballesteros}}, \bibinfo {author} {\bibfnamefont {A.}~\bibnamefont {Blasco}},\ and\ \bibinfo {author} {\bibfnamefont {I.}~\bibnamefont {Gutierrez-Sagredo}},\ }\bibfield  {title} {\enquote {\bibinfo {title} {Hamiltonian structure of compartmental epidemiological models},}\ }\href {https://doi.org/10.1016/j.physd.2020.132656} {\bibfield  {journal} {\bibinfo  {journal} {Phys. D: Nonlinear Phenom.}\ }\textbf {\bibinfo {volume} {413}} (\bibinfo {year} {2020}),\ 10.1016/j.physd.2020.132656}\BibitemShut {NoStop}%
\bibitem [{\citenamefont {O'Keeffe}, \citenamefont {Ceron},\ and\ \citenamefont {Petersen}(2022)}]{okeefe}%
  \BibitemOpen
  \bibfield  {author} {\bibinfo {author} {\bibfnamefont {K.}~\bibnamefont {O'Keeffe}}, \bibinfo {author} {\bibfnamefont {S.}~\bibnamefont {Ceron}},\ and\ \bibinfo {author} {\bibfnamefont {K.}~\bibnamefont {Petersen}},\ }\bibfield  {title} {\enquote {\bibinfo {title} {Collective behavior of swarmalators on a ring},}\ }\href@noop {} {\bibfield  {journal} {\bibinfo  {journal} {Phys. Rev. E}\ }\textbf {\bibinfo {volume} {105}} (\bibinfo {year} {2022})}\BibitemShut {NoStop}%
\bibitem [{\citenamefont {Pedergnana}\ and\ \citenamefont {Noiray}(2023)}]{pedergnana2023superradiant}%
  \BibitemOpen
  \bibfield  {author} {\bibinfo {author} {\bibfnamefont {T.}~\bibnamefont {Pedergnana}}\ and\ \bibinfo {author} {\bibfnamefont {N.}~\bibnamefont {Noiray}},\ }\bibfield  {title} {\enquote {\bibinfo {title} {Superradiant scattering by a limit cycle},}\ }\href {https://doi.org/10.1103/PhysRevApplied.20.034068} {\bibfield  {journal} {\bibinfo  {journal} {Phys. Rev. Applied}\ }\textbf {\bibinfo {volume} {20}},\ \bibinfo {pages} {034068} (\bibinfo {year} {2023})}\BibitemShut {NoStop}%
\bibitem [{\citenamefont {Sanders}, \citenamefont {Verhulst},\ and\ \citenamefont {Murdock}(2007)}]{sanders2007averaging}%
  \BibitemOpen
  \bibfield  {author} {\bibinfo {author} {\bibfnamefont {J.~A.}\ \bibnamefont {Sanders}}, \bibinfo {author} {\bibfnamefont {F.}~\bibnamefont {Verhulst}},\ and\ \bibinfo {author} {\bibfnamefont {J.}~\bibnamefont {Murdock}},\ }\href@noop {} {\emph {\bibinfo {title} {Averaging methods in nonlinear dynamical systems}}},\ Vol.~\bibinfo {volume} {59}\ (\bibinfo  {publisher} {Springer},\ \bibinfo {address} {Berlin, Heidelberg},\ \bibinfo {year} {2007})\BibitemShut {NoStop}%
\bibitem [{\citenamefont {Roberts}\ and\ \citenamefont {Spanos}(1986)}]{Roberts1986111}%
  \BibitemOpen
  \bibfield  {author} {\bibinfo {author} {\bibfnamefont {J.}~\bibnamefont {Roberts}}\ and\ \bibinfo {author} {\bibfnamefont {P.}~\bibnamefont {Spanos}},\ }\bibfield  {title} {\enquote {\bibinfo {title} {Stochastic averaging: An approximate method of solving random vibration problems},}\ }\href@noop {} {\bibfield  {journal} {\bibinfo  {journal} {Int. J. Non-Linear Mech.}\ }\textbf {\bibinfo {volume} {21}},\ \bibinfo {pages} {111--134} (\bibinfo {year} {1986})}\BibitemShut {NoStop}%
\bibitem [{\citenamefont {Noiray}(2017)}]{Noiray2017}%
  \BibitemOpen
  \bibfield  {author} {\bibinfo {author} {\bibfnamefont {N.}~\bibnamefont {Noiray}},\ }\bibfield  {title} {\enquote {\bibinfo {title} {Linear growth rate estimation from dynamics and statistics of acoustic signal envelope in turbulent combustors},}\ }\href {https://doi.org/10.1115/1.4034601} {\bibfield  {journal} {\bibinfo  {journal} {J. Eng. Gas Turbines Power}\ }\textbf {\bibinfo {volume} {139}} (\bibinfo {year} {2017}),\ 10.1115/1.4034601}\BibitemShut {NoStop}%
\bibitem [{\citenamefont {Balanov}\ \emph {et~al.}(2009)\citenamefont {Balanov}, \citenamefont {Janson}, \citenamefont {Postnov},\ and\ \citenamefont {Sosnovtseva}}]{balanov2009simple}%
  \BibitemOpen
  \bibfield  {author} {\bibinfo {author} {\bibfnamefont {A.}~\bibnamefont {Balanov}}, \bibinfo {author} {\bibfnamefont {N.}~\bibnamefont {Janson}}, \bibinfo {author} {\bibfnamefont {D.}~\bibnamefont {Postnov}},\ and\ \bibinfo {author} {\bibfnamefont {O.}~\bibnamefont {Sosnovtseva}},\ }\href@noop {} {\emph {\bibinfo {title} {From simple to complex}}}\ (\bibinfo  {publisher} {Springer},\ \bibinfo {address} {Berlin, Heidelberg},\ \bibinfo {year} {2009})\BibitemShut {NoStop}%
\bibitem [{\citenamefont {Van~der Pol}(1926)}]{van1926lxxxviii}%
  \BibitemOpen
  \bibfield  {author} {\bibinfo {author} {\bibfnamefont {B.}~\bibnamefont {Van~der Pol}},\ }\bibfield  {title} {\enquote {\bibinfo {title} {{LXXXVIII}. on “relaxation-oscillations”},}\ }\href@noop {} {\bibfield  {journal} {\bibinfo  {journal} {London Edinburgh Philos. Mag. J. Sci.}\ }\textbf {\bibinfo {volume} {2}},\ \bibinfo {pages} {978--992} (\bibinfo {year} {1926})}\BibitemShut {NoStop}%
\bibitem [{\citenamefont {Duffing}(1918)}]{duffing1918forced}%
  \BibitemOpen
  \bibfield  {author} {\bibinfo {author} {\bibfnamefont {G.}~\bibnamefont {Duffing}},\ }\href@noop {} {\emph {\bibinfo {title} {Erzwungene {S}chwingungen bei veränderlicher {E}igenfrequenz und ihre technische {B}edeutung}}},\ \bibinfo {number} {41-42}\ (\bibinfo  {publisher} {Vieweg},\ \bibinfo {year} {1918})\ pp.\ \bibinfo {pages} {1--134}\BibitemShut {NoStop}%
\bibitem [{\citenamefont {Songling}(1981)}]{Songling1981301}%
  \BibitemOpen
  \bibfield  {author} {\bibinfo {author} {\bibfnamefont {S.}~\bibnamefont {Songling}},\ }\bibfield  {title} {\enquote {\bibinfo {title} {A method of constructing cycles without contact around a weak focus},}\ }\href {https://doi.org/10.1016/0022-0396(81)90039-5} {\bibfield  {journal} {\bibinfo  {journal} {J. Differ. Equ.}\ }\textbf {\bibinfo {volume} {41}},\ \bibinfo {pages} {301 – 312} (\bibinfo {year} {1981})}\BibitemShut {NoStop}%
\bibitem [{\citenamefont {Kuznetsov}, \citenamefont {Kuznetsova},\ and\ \citenamefont {Leonov}(2013)}]{Kuznetsov201329}%
  \BibitemOpen
  \bibfield  {author} {\bibinfo {author} {\bibfnamefont {N.}~\bibnamefont {Kuznetsov}}, \bibinfo {author} {\bibfnamefont {O.}~\bibnamefont {Kuznetsova}},\ and\ \bibinfo {author} {\bibfnamefont {G.}~\bibnamefont {Leonov}},\ }\bibfield  {title} {\enquote {\bibinfo {title} {Visualization of four normal size limit cycles in two-dimensional polynomial quadratic system},}\ }\href {https://doi.org/10.1007/s12591-012-0118-6} {\bibfield  {journal} {\bibinfo  {journal} {Differ. Equ. Dyn. Syst.}\ }\textbf {\bibinfo {volume} {21}},\ \bibinfo {pages} {29 – 34} (\bibinfo {year} {2013})}\BibitemShut {NoStop}%
\bibitem [{\citenamefont {Dirac}(1981)}]{dirac1981principles}%
  \BibitemOpen
  \bibfield  {author} {\bibinfo {author} {\bibfnamefont {P.~A.~M.}\ \bibnamefont {Dirac}},\ }\href@noop {} {\emph {\bibinfo {title} {The principles of quantum mechanics}}},\ \bibinfo {number} {27}\ (\bibinfo  {publisher} {Oxford university press},\ \bibinfo {year} {1981})\BibitemShut {NoStop}%
\bibitem [{\citenamefont {Moiseyev}(2011)}]{moiseyev2011non}%
  \BibitemOpen
  \bibfield  {author} {\bibinfo {author} {\bibfnamefont {N.}~\bibnamefont {Moiseyev}},\ }\href@noop {} {\emph {\bibinfo {title} {Non-Hermitian quantum mechanics}}}\ (\bibinfo  {publisher} {Cambridge University Press},\ \bibinfo {year} {2011})\BibitemShut {NoStop}%
\bibitem [{\citenamefont {Deguchi}\ and\ \citenamefont {Fujiwara}(2020)}]{Deguchi2020}%
  \BibitemOpen
  \bibfield  {author} {\bibinfo {author} {\bibfnamefont {S.}~\bibnamefont {Deguchi}}\ and\ \bibinfo {author} {\bibfnamefont {Y.}~\bibnamefont {Fujiwara}},\ }\bibfield  {title} {\enquote {\bibinfo {title} {Quantization of the damped harmonic oscillator based on a modified {B}ateman {L}agrangian},}\ }\href {https://doi.org/10.1103/PhysRevA.101.022105} {\bibfield  {journal} {\bibinfo  {journal} {Phys. Rev. A}\ }\textbf {\bibinfo {volume} {101}} (\bibinfo {year} {2020}),\ 10.1103/PhysRevA.101.022105}\BibitemShut {NoStop}%
\bibitem [{\citenamefont {Blacker}\ and\ \citenamefont {Tilbrook}(2021)}]{Blacker2021}%
  \BibitemOpen
  \bibfield  {author} {\bibinfo {author} {\bibfnamefont {M.~J.}\ \bibnamefont {Blacker}}\ and\ \bibinfo {author} {\bibfnamefont {D.~L.}\ \bibnamefont {Tilbrook}},\ }\bibfield  {title} {\enquote {\bibinfo {title} {Alternative approach to the quantization of the damped harmonic oscillator},}\ }\href {https://doi.org/10.1103/PhysRevA.104.032211} {\bibfield  {journal} {\bibinfo  {journal} {Phys. Rev. A}\ }\textbf {\bibinfo {volume} {104}} (\bibinfo {year} {2021}),\ 10.1103/PhysRevA.104.032211}\BibitemShut {NoStop}%
\bibitem [{\citenamefont {Eichler}\ and\ \citenamefont {Zilberberg}(2023)}]{eichler2023classical}%
  \BibitemOpen
  \bibfield  {author} {\bibinfo {author} {\bibfnamefont {A.}~\bibnamefont {Eichler}}\ and\ \bibinfo {author} {\bibfnamefont {O.}~\bibnamefont {Zilberberg}},\ }\href@noop {} {\emph {\bibinfo {title} {Classical and Quantum Parametric Phenomena}}}\ (\bibinfo  {publisher} {Oxford University Press},\ \bibinfo {year} {2023})\BibitemShut {NoStop}%
\bibitem [{\citenamefont {Esen}, \citenamefont {Grmela},\ and\ \citenamefont {Pavelka}(2022)}]{Esen2022}%
  \BibitemOpen
  \bibfield  {author} {\bibinfo {author} {\bibfnamefont {O.}~\bibnamefont {Esen}}, \bibinfo {author} {\bibfnamefont {M.}~\bibnamefont {Grmela}},\ and\ \bibinfo {author} {\bibfnamefont {M.}~\bibnamefont {Pavelka}},\ }\bibfield  {title} {\enquote {\bibinfo {title} {On the role of geometry in statistical mechanics and thermodynamics. {I}. {G}eometric perspective},}\ }\href {https://doi.org/10.1063/5.0099923} {\bibfield  {journal} {\bibinfo  {journal} {J. Math. Phys.}\ }\textbf {\bibinfo {volume} {63}} (\bibinfo {year} {2022}),\ 10.1063/5.0099923}\BibitemShut {NoStop}%
\bibitem [{\citenamefont {Waalkens}, \citenamefont {Schubert},\ and\ \citenamefont {Wiggins}(2008)}]{Waalkens2008R1}%
  \BibitemOpen
  \bibfield  {author} {\bibinfo {author} {\bibfnamefont {H.}~\bibnamefont {Waalkens}}, \bibinfo {author} {\bibfnamefont {R.}~\bibnamefont {Schubert}},\ and\ \bibinfo {author} {\bibfnamefont {S.}~\bibnamefont {Wiggins}},\ }\bibfield  {title} {\enquote {\bibinfo {title} {Wigner's dynamical transition state theory in phase space: Classical and quantum},}\ }\href {https://doi.org/10.1088/0951-7715/21/1/R01} {\bibfield  {journal} {\bibinfo  {journal} {Nonlinearity}\ }\textbf {\bibinfo {volume} {21}},\ \bibinfo {pages} {R1–R118} (\bibinfo {year} {2008})}\BibitemShut {NoStop}%
\bibitem [{\citenamefont {Jahn}\ \emph {et~al.}(2020)\citenamefont {Jahn}, \citenamefont {Stender}, \citenamefont {Tatzko}, \citenamefont {Hoffmann}, \citenamefont {Grolet},\ and\ \citenamefont {Wallaschek}}]{Jahn2020}%
  \BibitemOpen
  \bibfield  {author} {\bibinfo {author} {\bibfnamefont {M.}~\bibnamefont {Jahn}}, \bibinfo {author} {\bibfnamefont {M.}~\bibnamefont {Stender}}, \bibinfo {author} {\bibfnamefont {S.}~\bibnamefont {Tatzko}}, \bibinfo {author} {\bibfnamefont {N.}~\bibnamefont {Hoffmann}}, \bibinfo {author} {\bibfnamefont {A.}~\bibnamefont {Grolet}},\ and\ \bibinfo {author} {\bibfnamefont {J.}~\bibnamefont {Wallaschek}},\ }\bibfield  {title} {\enquote {\bibinfo {title} {The extended periodic motion concept for fast limit cycle detection of self-excited systems},}\ }\href {https://doi.org/10.1016/j.compstruc.2019.106139} {\bibfield  {journal} {\bibinfo  {journal} {Comput. Struct.}\ }\textbf {\bibinfo {volume} {227}} (\bibinfo {year} {2020}),\ 10.1016/j.compstruc.2019.106139}\BibitemShut {NoStop}%
\end{thebibliography}%

\end{document}